%% file: Main.tex
\documentclass[12pt,fleqn,reqno]{article}
\usepackage{amsmath}
\usepackage{amssymb}
\usepackage{amsfonts}
\usepackage{amsthm}
\usepackage{mathtools}
\mathtoolsset{
above-intertext-sep = -1ex
}

\usepackage[utf8]{inputenc}
\usepackage[T1]{fontenc}

\usepackage[sans]{dsfont}

\usepackage[left=2.5cm, right=2.5cm, top=2.5cm, bottom=2.5cm]{geometry}
\usepackage[%
bb=fourier,
scr=euler,
]
{mathalfa}

\usepackage[font=small,labelfont=bf]{caption}
\usepackage{subfigure}
\usepackage{graphicx}
\graphicspath{{Figures/}} 
\usepackage{acronym}
\usepackage{latexsym}
\usepackage{paralist}
\usepackage{wasysym}
\usepackage{xspace}
\usepackage{framed}
\usepackage{palatino,pxfonts}
\usepackage{authblk}
\usepackage[numbers,sort&compress]{natbib}

\usepackage[dvipsnames,svgnames]{xcolor}
\colorlet{MyBlue}{DodgerBlue!75!Black}
\colorlet{MyGreen}{DarkGreen!95!Black}

\usepackage{hyperref}
\hypersetup{
colorlinks=true,
linktocpage=true,
pdfstartview=FitH,
breaklinks=true,
pdfpagemode=UseNone,
pageanchor=true,
pdfpagemode=UseOutlines,
plainpages=false,
bookmarksnumbered,
bookmarksopen=false,
bookmarksopenlevel=1,
hypertexnames=true,
pdfhighlight=/O,
urlcolor=MyBlue!60!black,linkcolor=MyBlue!70!black,citecolor=DarkGreen!70!black, % <--- for screen
pdftitle={},
pdfauthor={},
pdfsubject={},
pdfkeywords={},
pdfcreator={pdfLaTeX},
pdfproducer={LaTeX with hyperref}
}
\numberwithin{equation}{section}  %numberwithin goes before cleverefs when using hyperref
\usepackage[sort&compress,capitalize,nameinlink]{cleveref}
\crefname{example}{Ex.}{Exs.}

\crefrangeformat{equation}{\upshape(#3#1#4)\textendash(#5#2#6)}

\newcommand{\dd}{\:d}

\newcommand{\eps}{\varepsilon}

\newcommand{\dif}{\dd}
\newcommand{\bigoh}{\mathcal{O}}

\DeclareMathOperator*{\argmin}{argmin}

\DeclareMathOperator{\bd}{bd}
\DeclareMathOperator{\cl}{cl}
\DeclareMathOperator{\epi}{epi}

\DeclareMathOperator{\diag}{diag}

\DeclareMathOperator{\dist}{dist}

\DeclareMathOperator{\dom}{dom}

\DeclareMathOperator{\Int}{ri}

\DeclareMathOperator{\grad}{grad}
\DeclareMathOperator{\image}{im}

\DeclareMathOperator{\tr}{tr}

\DeclareMathOperator{\Id}{Id}

\newcommand{\const}{\mathtt{c}}

\newcommand{\ca}{\mathtt{a}}
\newcommand{\cb}{\mathtt{b}}
\DeclareMathOperator{\lev}{lev}
\newcommand{\bK}{\mathbf{K}}
\newcommand{\bN}{\mathbf{N}}
\renewcommand{\iff}{\Leftrightarrow}

\renewcommand{\emptyset}{\varnothing}

\newcommand{\scrA}{\mathcal{A}}

\newcommand{\scrC}{\mathcal{C}}

\newcommand{\scrE}{\mathcal{E}}
\newcommand{\scrF}{\mathcal{F}}
\newcommand{\scrG}{\mathcal{G}}
\newcommand{\scrH}{\mathcal{H}}

\newcommand{\scrL}{\mathcal{L}}
\newcommand{\scrM}{\mathcal{M}}

\newcommand{\scrS}{\mathcal{S}}

\newcommand{\scrW}{\mathcal{W}}
\newcommand{\scrX}{\mathcal{X}}

\newcommand{\Rn}{\R^n}

\newcommand{\R}{\mathbb{R}}

\newcommand{\N}{\mathbb{N}}
\DeclareMathOperator{\NC}{\mathsf{NC}}
\DeclareMathOperator{\TC}{\mathsf{TC}}

\newcommand{\bC}{{\mathbf{C}}}

\newcommand{\metric}{\mathsf{d}}

\usepackage{algorithm2e}
\usepackage{algpseudocode}								% for algorithm macros

\DeclareMathOperator{\HBA}{HBA}
\DeclareMathOperator{\AHBA}{AHBA}
\theoremstyle{plain}
\newtheorem{theorem}{Theorem}
\newtheorem{corollary}[theorem]{Corollary}
\newtheorem*{corollary*}{Corollary}
\newtheorem{lemma}[theorem]{Lemma}
\newtheorem{proposition}[theorem]{Proposition}

\theoremstyle{definition}
\newtheorem{definition}[theorem]{Definition}
\newtheorem*{definition*}{Definition}
\newtheorem{assumption}{Assumption}
\newenvironment{Proof}[1][Proof]{\begin{proof}[#1]}{\end{proof}}

\theoremstyle{remark}
\newtheorem{remark}{Remark}
\newtheorem*{remark*}{Remark}
\newtheorem*{notation*}{Notational remark}

\numberwithin{theorem}{section}
\numberwithin{remark}{section}
\numberwithin{example}{section}

\DeclarePairedDelimiter{\abs}{\lvert}{\rvert}

\DeclarePairedDelimiter{\inner}{\langle}{\rangle}
\DeclarePairedDelimiter{\norm}{\lVert}{\rVert}

\newacro{HBA}[HBA]{Hessian-barrier algorithm}
\newacro{AHBA}[AHBA]{Adaptive-Hessian-barrier algorithm}
\newacroplural{NE}[NE]{Nash equilibria}
\newacro{VI}{variational inequality}
\newacroplural{VI}[VIs]{variational inequalities}
\newacro{iid}[i.i.d.]{independent and identically distributed}
\newacro{FOM}{First-order method}

\title{Generalized Self-concordant Hessian-barrier algorithms}
\date{\today}

\author[1]{\small Pavel Dvurechensky}
\author[2]{\small Mathias Staudigl}
\author[3]{\small Cesar A. Uribe}

\affil[1]{\footnotesize Weierstrass Institute for Applied Analysis and Stochastics, Mohrenstr. 39, 10117 Berlin, Germany\\
(\href{mailto:Pavel.Dvurechensky@wias-berlin.de}{Pavel.Dvurechensky@wias-berlin.de})}
\affil[2]{\footnotesize Maastricht University, Department of Quantitative Economics, P.O. Box 616, NL\textendash 6200 MD Maastricht, The Netherlands\\
(\href{mailto:m.staudigl@maastrichtuniversity.nl}{m.staudigl@maastrichtuniversity.nl})}
\affil[3]{\footnotesize Laboratory for Information and Decision Systems (LIDS)\\ 
Institute for Data, Systems, and Society (IDSS), Massachusetts Institute of Technology, Cambridge MA\\
(\href{mailto:cauribe@mit.edu}{cauribe@mit.edu})}
\begin{document}
\maketitle
%%% ABSTRACT
%----------------------------------------------------------------------
\begin{abstract}
\input{abstract}
\end{abstract}

%*************************************************************st
%*****    BODY TEXT
%*************************************************************
\renewcommand{\sharp}{\gamma}
\acresetall
\allowdisplaybreaks

%----------------------------------------------------------------------
%%% INTRODUCTION
%----------------------------------------------------------------------
\section{Introduction}
\label{sec:introduction}
\input{Introduction}

%----------------------------------------------------------------------
%%% PRELIMS
%----------------------------------------------------------------------
\section{Setup and preliminaries}
\label{sec:prelims}
\input{prelims}
\section{The Hessian-barrier method}
\label{sec:method}
%\subsection{Defining the search directions}
\label{sec:general}
\input{HBA_general}
\subsection{The Hessian-Barrier potential reduction algorithm}
\label{sec:PotentialReduction}
\input{PotentialReduction}
%----------------------------------------------------------------------
%% Algorithm
%----------------------------------------------------------------------
\section{The Hessian-barrier algorithm}
\label{sec:algo}
\input{algo}

%----------------------------------------------------------------------
%%% Convergence
%----------------------------------------------------------------------
\section{Complexity analysis of HBA}
\label{sec:complexity}
\input{convergence}

%----------------------------------------------------------------------
%%% Numerics
%----------------------------------------------------------------------
\section{Numerical Results}
\label{sec:numerics}
\subsection{Statistical learning with non-convex regularization}
\input{SCAD}

\subsection{$L^{p}$-minimization}
\label{sec:Lp}
\input{Lp}
%\subsection{Poisson linear inverse problems}
%\label{sec:Poisson}
%\input{./Applications/Poisson}
%----------------------------------------------------------------------
%%% CONCLUSIONS
%----------------------------------------------------------------------
\section{Conclusion}
\label{sec:conclusion}
\input{conclusion}

%----------------------------------------------------------------------
%%% CONCLUSIONS
%----------------------------------------------------------------------
\section{Conclusion}
\label{sec:conclusion}
\input{conclusion}
%*************************************************************
%*****    ACKNOWLEDGMENTS
%*************************************************************
\section*{Acknowledgments.}
M. Staudigl would like to thank Panayotis Mertikopoulos for years of fruitful collaborations, and extensive feedback on this paper. The research of M. Staudigl has been supported by the COST Action CA16228 "European Network for Game Theory". The work of C.A. Uribe was partially supported by Yahoo! Research Faculty Engagement Program. The work by P. Dvurechensky was supported by RFBR grants 18-31-20005 mol\_a\_ved and 18-29-03071\_mk.

%*************************************************************
%*****    APPENDICES
%*************************************************************
\begin{appendix}
%----------------------------------------------------------------------
%%% APP1:  
%----------------------------------------------------------------------
\section{Proof of Proposition \ref{prop:step}}
\label{app:Step}
\input{Appendix_Step}

\section{Proof of Theorem \ref{th:gradient}}
\label{app:stationary}
\input{appendix_stationary}
\end{appendix}

%*************************************************************
%*****    BIBLIOGRAPHY
%*************************************************************
\bibliographystyle{plainnat}
\bibliography{mybib}
\end{document}

%% file: abstract.tex
%%% Abstract
%----------------------------------------------------------------------
% !TEX root = ./Main.tex
%

Many problems in statistical learning, imaging, and computer vision involve the optimization of a non-convex objective function with singularities at the boundary of the feasible set. For such challenging instances, we develop a new interior-point technique building on the Hessian-barrier algorithm recently introduced in Bomze, Mertikopoulos, Schachinger and Staudigl, [SIAM J. Opt. 2019 29(3), pp. 2100-2127], where the Riemannian metric is induced by a generalized self-concordant function. This class of functions is sufficiently general to include most of the commonly used barrier functions in the literature of interior point methods. We prove global convergence to an approximate stationary point of the method, and in cases where the feasible set admits an easily computable self-concordant barrier, we verify worst-case optimal iteration complexity of the method. Applications in non-convex statistical estimation and $L^{p}$-minimization are discussed to given the efficiency of the method.

%% file: Introduction.tex
%%% Introduction
%----------------------------------------------------------------------
% !TEX root = ./Main.tex
%

In this paper, we consider the following constrained minimization problem, which has plenty of applications in diverse disciplines, including machine learning, signal processing, statistics, and operations research 
\begin{equation}\label{eq:P}\tag{P}
f^{\ast}=\min\{f(x): x\in \bar{\scrC},Ax=b\}
\end{equation}
Here $\bar{\scrC}\subset\Rn$ is a nonempty, closed and convex set, and $f$ is a (possibly) non-smooth, non-convex function from $\Rn\to\R\cup\{+\infty\}$. A special case of \eqref{eq:P} are regularized statistical estimation problems where the aim is to find a parameter vector $x\in\Rn$ in order to minimize a composite objective of the form 
\begin{equation}\label{eq:model}
f(x):=f_{0}(x)+f_{1}(c(x)).
\end{equation}
In such applications the function $f_{0}:\Rn\to\R$ is a continuous data fidelity term, $f_{1}:\R^{d}\to(-\infty,\infty]$ is a regularizer and $c:\Rn\to\R^{d}$ is some link function mapping the parameters to a usually lower dimensional subspace. Typical formulations of such problems can be given in the form of 
\begin{equation}\label{eq:L2LP}
f(x)=\frac{1}{2}\norm{s-Wx}^{2}_{2}+\sum_{i=1}^{m}\varphi_{i}(\norm{D_{i}x}^{p})
\end{equation}
in which $f_{0}(x)=\frac{1}{2}\norm{s-Wx}_{2}^{2}, f_{1}(y)=\sum_{i=1}^{m}\varphi_{i}(y_{i}):\R^{q}\to[-\infty,\infty)$, $\varphi_{i}:\R^{d_{i}}\to\R$ is continuously differentiable, and $c(x)=(\norm{D_{1}x}^{p},\ldots,\norm{D_{m}x}^{p})^{\top}$ for some $p\in(0,1),D_{i}\in\R^{d_{i}\times n},d=\sum_{i=1}^{m}d_{i}$. In fact, the formulation \eqref{eq:L2LP} includes many well-known problems in statistics: fused lasso \cite{TibRosSauZhu05}, grouping pursuit \cite{SheHua10}, etc. From an optimization perspectives, these regularized estimation problems are challenging since they are non-convex and not globally Lipschitz continuous and thus belong to the class of NP-hard problems. Even worse, they may even fail to be differentiable. As an example, the $L^{p}$ regularization problem with link function $c(x)=(\abs{x_{1}}^{p},\ldots,\abs{x_{n}}^{p})$ for $0<p<1$ and $f_{1}(y)=\sum_{i=1}^{n}y_{i}$ fails to be even directionally differentiable when $x_{i}=0$ for some $i=1,2,\ldots,n$.

A common tenet of all recent applications of \eqref{eq:P} is that the problem involves a huge number of variables. This makes the application of classical interior-point solvers infeasible. Instead, first-order methods (FOMs) with cheap per-iteration implementation costs are the method of choice \cite{Teb18}. The most impactful success stories of FOMs have been achieved under the quite demanding assumption that the objective function is convex and smooth, and the feasible set admits a proximal-friendly formulation. Indeed, if proximal-based projection operators onto the feasible set $\scrX:=\{x\in\bar{\scrC}: Ax=b\}$ are easy to evaluate, black-box based FOMs such as mirror descent, projected subgradient, and conditional gradient methods can be tuned to successfully solve \eqref{eq:P} up to $\eps$-accuracy. How to handle non-convex objective functions with first-order methods is still a challenging problem receiving a lot of interest from various different perspectives, in particular in statistical and deep learning. Beside this implicit assumption in all projection-based FOMs, another fundamental assumption in all these methods is the availability of a \emph{descent lemma} \cite[Lem.1.2.3]{NesConvex}. Sufficient for such an a-priori estimate is that the objective function's gradient is a Lipschitz continuous function. Already the above mentioned application to statistical estimation shows that this Lipschitz-smoothness assumption is to demanding to cover such important applications. Only recently, the path-breaking paper \cite{BauBolTeb16} resolved this problem by introducing the concept of \emph{relative smoothness} as a surrogate for the demanding Lipschitz gradient assumption (see also \cite{LuFreNes18} for elaborations and applications). Based on relative-smoothness, they derive a new descent lemma where the usual quadratic approximation is replaced by a non-Euclidean proximity measure, which captures the objective function and the geometry of the underlying domain all at once. The corresponding proximal-based subgradient algorithm comes with global convergence guarantees and complexity estimates for \emph{convex composite} models. This beautiful, and practically relevant, approach has been recently extended to a non-convex composite model in \cite{BolSabTebVai18}, where new complexity estimates are derived as well.

\subsection{Our Approach}
This work is concerned with a different approach to tackle non-convex optimization problems avoiding knowledge of a global Lipschitz constant. Our work is inspired by the recent Riemannian gradient methods developed in \cite{HBA-linear}, where a rather large class of \emph{Hessian barrier algorithms} (HBA) has been constructed as numerical schemes for solving Lipschitz-smooth, non-convex optimization problems over the polyhedron $\{x\in\Rn: x\geq 0,Ax=b\}$. The construction of HBA is motivated by looking at an \emph{explicit} numerical discretization of a continuous-time dynamical systems introduced by \cite{ABB04,BolTeb03} as a theoretical method for solving \emph{convex} linearly constrained smooth optimization problems. However, these schemes remained at a conceptual level, and the usefulness of these dynamical systems for effectively solving constrained optimization of the form \eqref{eq:P} remained completely unanswered. $\HBA$ laid the foundations to an algorithmic analysis for these dynamical systems, and investigated their efficiency when solving linearly constrained and smooth non-convex optimization problems. $\HBA$ first identifies the feasible set as a Riemannian manifold with a metric induced by the Hessian matrix of a $\bC^{2}$ barrier-like function $h$ (a \emph{barrier-generating kernel}). Once the geometry has been defined, a step-size strategy is designed ensuring feasibility and a sufficient decrease of the objective. The analysis in \cite{HBA-linear} relied heavily on the classical Lipschitz-descent lemma \cite{NesConvex}, and involved an Armijo line-search procedure. It has been shown that this approach generalizes many classical interior point techniques like affine scaling \cite{BayLag89}, and also contains Lotka-volterra systems as a special case \cite{TseBomSch11}. A complexity analysis was achieved in the case where the objective function is quadratic, and it has been shown that a proper choice of the Riemannian metric affects the complexity of the method \cite[Thm 5.1]{HBA-linear}. In this paper we significantly extend the results obtained in \cite{HBA-linear} by constructing a new first-order interior point method for solving problem \eqref{eq:P} under very mild assumptions on the data. In order to explain the approach described in this paper, let us go back to the classical way of solving problem \eqref{eq:P} when $f$ is \emph{convex}. The most famous algorithm for solving such problems are interior point methods, which solve conic optimization problems in polynomial time \cite{NesNem94}. The key structure exploited in conventional IPMs is the existence of a \emph{self-concordant barrier} (SCB) for the set constraint $\bar{\scrC}$. In such cases one considers the potential function 
\[
F_{\mu}(x)=f(x)+\mu h(x),
\]
where $\mu>0$ is a penalty parameter and $h$ is a barrier function over the set $\bar{\scrC}$. By fixing a sequence of barrier parameters $(\mu_{k})_{k\geq 0}$ with $\mu_{k}\downarrow 0$ and solving the sequence of minimization problems $\min_{x}F_{\mu_{k}}(x)$ along this sequence generates the \emph{analytic central path} $\{x^{\ast}_{\mu}:\mu>0\}$ as it converges to the solution of the actual problem of interest \eqref{eq:P}. For proving convergence of the analytic central path, SCBs are the key tool to prove polynomial-solvability of the barrier problem by sequentially using Newton's method.\footnote{Recently, the path-following method was extended for self-concordant functions which are not self-concordant barriers in \cite{DvuNes18}.} The new Hessian-barrier method we propose in this paper follows similar ideas. We first embed the original optimization problem \eqref{eq:P} into a potential-reduction scheme involving the potential function $F_{\mu}$. However, instead of classical self-concordance theory our analysis works easily on a much broader class of \emph{generalized self-concordant functions} (GSC), which have recently been introduced in \cite{SunTran18}. As we show in this paper, GSC functions provide a very attractive class of barrier-generating kernels as their Hessian matrix induce a Riemannian metric under which a full-fledged complexity analysis can be performed. 

\subsection{Available complexity results}
We now review results on complexity analysis of non-convex optimization problems. Since the Hessian-barrier method uses only first-order information about $f$, we restrict our review to first-order methods as well. In unconstrained non-convex optimization the usual criticality measure is the norm of the gradient of the objective function. Hence, oracle complexity of a given algorithm refers to the number of oracle queries until $\norm{\nabla f(x^{k})}\leq\eps$ for a targeted tolerance level $\eps>0$. An in-depth survey of known complexity bounds can be found in \cite{CarDucHinSid19,CarDucHinSid19b}. We can give no justice to the huge literature on complexity results for first-order methods but provide below a partial survey of known complexity estimates in order to put our results into perspective.  

\paragraph{\emph{Smooth, non-convex}} For quadratic programming problems with linear constraints, the authors in \cite{Ye98} proved that an $\eps$-KKT point is computed in $\bigoh(\eps^{-1}\log(\eps^{-1}))$ iterations. A recent manifestation of the effect of Riemannian geometry on the convergence to $\eps$-KKT points can be found in \cite{HBA-linear}. For general unconstrained nonconvex optimization, it was shown in \cite{NesConvex} that a steepest descent with line search method can find an $\eps$-stationary point in $\bigoh(\eps^{-2})$ iterations. An accelerated method with the same guarantee can be found in \cite{guminov2019accelerated}. 
The same worst case complexity result holds for trust-region methods \cite{GraSarToi08}. The results for gradient method were generalized to the case of simple projection-friendly constraints and H\"older derivatives in \cite{ghadimi2019generalized}. 
Accelerated methods with complexity $\widetilde{\bigoh}(\eps^{-7/4})$ under additional assumption of Lipschitz second derivative are proposed in \cite{agarwal2017finding,carmon2017convex,carmon2018accelerated,jin2018accelerated}.

\paragraph{\emph{Lipschitz continuous, nonconvex}} Cartis, Gould and Toint \cite{CarGouToi11} estimated the worst-case complexity of a first-order trust-region or quadratic regularization method for solving unconstrained, non-convex minimization problems of the form \eqref{eq:model}, where $f_{1}:\R^{m}\to\R$ is convex but may be nonsmooth and $c:\Rn\to\R^{d}$ is continuously differentiable. Their method takes at most $\bigoh(\eps^{-2})$ iterations to reduce the size of a suitably defined first-order criticality measure below $\eps$.

\paragraph{\emph{Non-Lipschitz, nonconvex}} The authors in \cite{GeJiaYe11} extended the complexity result of \cite{Ye98} to the $L^{p}$-minimization problem over a polytop. They showed that finding and $\eps$-scaled stationary point or global minimizer requires at most $\bigoh(\eps^{-1}\log(\eps^{-1}))$ iterations. For general linear constrained non-convex minimization problems \cite{HaeLiuYe18} obtained an iteration complexity of $\bigoh(\eps^{-2})$ to reach an $\eps$-KKT point for optimization problems whose feasible set is defined by linear equality and non-negativity constraints. In the case of $L^{p}$ minimization for $p\in(0,1)$ over box constraints \cite{BiaCheYe15} develop a first-order interior point method yielding $\bigoh(\eps^{-2})$ worst-case iteration complexity in order to return and $\eps$-scaled stationary point.

\subsection{Our Contribution}
In relation to the above summarized literature, we provide here an easy-to-implement first-order method for non-convex non-smooth optimization problems, without requiring knowledge on the Lipschitz constant of the objective function. Specifically, the main contributions of this paper are summarized as follows:
\begin{enumerate}
\item We provide a new first-order interior point method based on the $\HBA$ method for non-convex and non-smooth optimization problems \eqref{eq:P} without Lipschitz smoothness conditions.
\item We are the first who provide a first-order interior point analysis based on GSC functions. 
\item We show how some model parameters can be made adaptive, making the method even more efficient in practice.
\item Our method comes with an explicit construction of optimal step-size policies and convergence guarantees.
\item We demonstrate optimal iteration complexity estimates on the order of $\bigoh(\eps^{-2})$ to reach a generalized stationary point, and connect this to classical $\eps$-KKT conditions in case where the barrier-generating kernel is a SC-B. This answers an open question raised in Remark 4.1. in \cite{BolTeb03}, since HBA can be seen as a descendent of the \emph{A-driven descent methods} defined in that paper. Also, in view of the partial literature survey given above, this rate is optimal.
%\item We test our method on two challenging non-convex formulati, and demonstrate the efficiency of $\HBA$. 
\end{enumerate}
The rest of this paper is organized as follows: Section \ref{sec:prelims} introduces the standing assumptions used in this paper, and introduces the class of generalized self-concordant functions. Section \ref{sec:method} defines conceptually the $\HBA$ method and introduces the optimal step-size policy associated with it. An adaptive variant of this base scheme is discussed in Section \ref{sec:adaptive}. Section \ref{sec:complexity} includes the main results in terms of convergence and complexity of the method. Section \ref{sec:numerics} contains numerical examples.

\paragraph{\emph{Notation}}
Given a $k$-times continuously differentiable function $f:\scrC\to\R$ and vectors $v_{1},\ldots,v_{k}\in\Rn$. For $x\in\Int(\scrC)$ and $1\leq j\leq k$, we define recursively 
\begin{align*}
D^{j}f(x)[v_{1},\ldots,v_{j}]:=\lim_{\eps\to 0}\frac{D^{j-1}f(x+\eps v_{j})[v_{1},\ldots,v_{j-1}]-D^{j-1}f(x)[v_{1},\ldots,v_{j-1}]}{\eps}.
\end{align*}
As a convention $D^{0}f(x)=f(x)$ and for $k=1$ we recover the directional derivative $f'(x;v)$. Given a positive semi-definite matrix $H\in\R^{n\times n}$, we define the norm $\norm{a}_{H}:=\sqrt{\inner{Ha,a}}$, and the dual norm $\norm{a}^{\ast}_{H}=\sup\{\inner{a,d}: \norm{d}_{H}=1\}$. If $H$ is invertible, the dual norm is given by $\norm{a}^{\ast}_{H}=\sqrt{\inner{H^{-1}a,a}}$. For a given $n\times n$ matrix $A$, let us define the operator norm $\abs{A}:=\sup\{\norm{Ax}:\norm{x}=1\}.$ Let $\scrC\subset\Rn$ be a convex set with closure $\bar{\scrC}$. Define the tangent cone $\TC_{\scrC}(x):=\cl\left[\R_{+}(\bar{C}-x)\right]$, and the corresponding polar cone $\NC_{\bar{\scrC}}(x)=\{p\in\Rn: \inner{p,y-x}\leq 0,\forall y\in\bar{C}\}$, which is called the normal cone.

%% file: prelims.tex
%%% prelims
%----------------------------------------------------------------------
% !TEX root = ./Main.tex
%

%In order to derive a family of first-order methods by following the HBA methodology, we need to accomplish two tasks: First, we need to find a good metric structure so that we can identify a computable step-size policy and thus be able to follow the computed search directions without violating primal feasibility. Second, we need to connect the asymptotic analysis of the HBA method with the underlying optimization problem, so that we can be sure to reach stationary points of problem \eqref{eq:Opt}. To meet this balance, we introduce below a class of metric-inducing kernels serving as good generators of our algorithmic scheme. 
%----------------------------------------------------------------------
%%% RIEMANNIAN
%----------------------------------------------------------------------
\subsection{Elements of Riemannian geometry}
\label{sec:Riemannian}

A key notion in our considerations is that of a \emph{Riemannian metric}, i.e. a position-dependent variant of the ordinary (Euclidean) scalar product between vectors \cite{Lee97}.
To define it, recall first that a \emph{scalar product} on $\Rn$ is a symmetric, positive-definite bilinear form $\inner{\cdot,\cdot}:\Rn\times \Rn\to \R$. This scalar product defines a norm in the usual way and it can be represented equivalently via its \emph{metric tensor}, that is, a symmetric, positive-definite matrix $H\in\R^{n\times n}$ with components
\begin{equation}
\label{eq:metric-components}
H_{ij}= \inner{e_{i},e_{j}}
\end{equation}
in the standard basis $\{e_i\}_{i=1}^{n}$ of $\Rn$. A \emph{Riemannian metric} on a nonempty open set $\scrC\subseteq\Rn$ is then defined to be a smooth assignment of scalar products $\inner{\cdot,\cdot}_{x}$ to each $x\in \scrC$ or, equivalently, a smooth field $H(x)$ of symmetric positive-definite matrices on $\scrC$.

Given a Riemannian metric on $\scrC$, the \emph{Riemannian gradient} of a smooth function $\phi:\scrC\to\R$ at $x\in \scrC$ is defined via the characterization
\begin{equation}
\label{eq:grad-char}
\inner{\grad\phi(x),z}_{x}= \phi'(x;z)	\quad	\text{for all $z\in\Rn$}.
\end{equation}
More concretely, by expressing everything in components, it is easy to see that
$\grad\phi(x)$ is given by the explicit expression 
\begin{equation}
\label{eq:grad-def}
\grad\phi(x)= [H(x)]^{-1} \nabla\phi(x).
\end{equation}

Bringing the above closer to our setting, let $\scrA_{0}\subseteq\Rn$ be a subspace of $\Rn$ and let $\scrA$ be an affine translate of $\scrA_{0}$ such that $\scrX^{\circ} \equiv C\cap\scrA$ is nonempty.
Then, viewing $\scrX^{\circ}$ as an open subset of $\scrA$, the \emph{gradient of $\phi$ restricted to $\scrX^{\circ}$} is defined as the unique vector $\grad_{\scrX^{\circ}} \phi(x) \equiv \grad \phi\vert_{\scrX^{\circ}}(x) \in\scrA_{0}$ such that
\begin{equation}
\label{eq:grad-res}
\grad_{\scrX^{\circ}}\phi(x)=\phi'(x,d)	\quad	\text{for all $d\in\scrA_{0}$}.
\end{equation}
Hence, specializing all this to the problem at hand, let $H(x)$ be a Riemannian metric on the open set $\scrC\subset\Rn$ and set
\[
\scrA_{0}:= \ker(A):= \{d\in\Rn: Ad=0\},\quad \scrA:=\{x\in\Rn: Ax=b\}.
\]
Then, a straightforward exercise in matrix algebra shows that the gradient of $f$ restricted to $\scrX^{\circ}=\scrC\cap\scrA_{0}$ can be written in closed form as
\begin{equation}
\label{eq:grad-res-coords}
\grad_{\scrX^{\circ}}f(x)= P_{x}[H(x)]^{-1} \nabla f(x)
\end{equation}
with 
\begin{equation}\label{eq:project}
P_{x}:= \Id - [H(x)]^{-1}A^{\top}(A[H(x)]^{-1}A^{\top})^{-1}A.
\end{equation}

%----------------------------------------------------------------------
%%% Prox
%----------------------------------------------------------------------
\subsection{Generalized self-concordant functions}
\label{sec:SC}
The Hessian-barrier method's main assumption is that the set constraint $\bar{\scrC}$ admits an easy-to-compute \emph{generalized self-concordant function}. In the next sections we describe the subclass of admissible functions on which the subsequent constructions build on. We begin with the notion of \emph{kernel generating distance}, by following the very general setup introduced in \cite{AusTeb06}. 
%%%%%%%%%%%%%%%%%%%%%%%%%%%%%%%%%%%%%
\begin{definition}[Kernel generating distance] 
\label{def:kgd}
Let $\scrC$ be a nonempty convex and open subset of $\Rn$. Associated with $\scrC$, a function $h:\scrC\to(-\infty,\infty]$ is called a \emph{kernel generating distance} if
\begin{itemize}
\item[(a)] $h$ is proper, lower semi-continuous and convex, with $\dom h\subset\bar{\scrC}$ and $\dom\partial h=\Int\dom h=\scrC$. 
\item[(b)] $h\in\bC^{1}(\scrC)$. 
%\item[(c)] $\bar{\scrC}\ni x\mapsto \norm{\nabla h(x)}^{\ast}_{x}$ is continuous. 
\end{itemize}
Denote the class of kernel generating distances by $\scrG(\scrC)$.
\end{definition}
%%%%%%%%%%%%%%%%%%%%%%%%%%%
%Properties (a) and (b) are quite standard in perceived proximal methods. Property (c) is not standard in this framework but needed in our analysis in order to align the long-run behavior of the algorithm with KKT conditions. We will spent some more words on the interpretation of this condition in Remark \ref{rem:dual}. 

%Given $h\in\scrG(\scrC)$, define the \emph{Bregman distance} $D_{h}:\dom(h)\times \scrC\to\R_{+}$ as 
%\begin{equation}\label{eq:Bregman}
%D_{h}(x,y)=h(x)-h(y)-\inner{\nabla h(y),x-y}
%\end{equation}
%By definition of the function $h\in\scrG(\scrC)$, we have 
%\begin{align*}
%D_{h}(x,y)\geq \frac{\kappa_{h}}{2}\norm{x-y}^{2}\qquad\forall (x,y)\in\dom(h)\times \scrC
%\end{align*}
%where $\kappa_{h}\geq 0$ is the strong convexity constant of the function $h$.
%%%%%%%%%%%%%%%%%%%%%%%%%%%%%%%%%%%%%%%%%%%%%%%%%%%%%%
%----------------------------------------------------------------------
%%% GSC
%----------------------------------------------------------------------
%\subsection{Generalized self-concordant functions}
%\label{sec:SC}
We now add the additional structure on the kernel generating distance $h\in\scrG(\scrC)$ we use in our algorithmic design. %The class of methods we develop in this paper follows the basic HBA methodology \cite{HBA-linear}. The key driver of this method is the Riemannian metric with which the affine manifold $\scrX$ is endowed. Following the the class of Riemannian dynamical systems \cite{ABB04,BT03,AttTeb04}, our aim here is develop a bona-fide discrete-time numerical algorithm containing all these dynamics as specific and \emph{explicit} time discretizations, and (b) provide for the first time a complexity analysis for the thus derived methods. We should point out that this answers an open problem raised in \cite{BT03}, where a large class of Barrier operators is constructed serving as a generator for solving smooth convex optimization problems. The Riemannian metric in which the Hessian-barrier method discussed in this paper works is derived from the class of \emph{generalized self-concordant functions} (GSC), as recently introduced in \cite{SunTran18}. Generalized self-concordance is a huge extension of the classical concept of self-concordant functions \cite{NesNemIPM94}, fundamental to the theory of interior-point methods. We recall the definition of GSC functions in this section, provide some examples, and discuss key estimates needed in this paper. 
%%%%%%%%%%%%%%%%%%%%%%%%%%%%%%%%%%%%%
\begin{definition}\citep{SunTran18}
Let $\phi\in\bC^{3}(\dom\phi;\R)$ be a closed convex function with $\dom\phi$ open. Given $\nu>0$ and $M>0$ some constants, we say that $\phi$ is $(M,\nu)$ is \emph{generalized self-concordant} (GSC) if 
\begin{equation}
\abs{\phi'''(t)}\leq M\phi''(t)^{\frac{\nu}{2}}\qquad\forall t\in\dom\phi. 
\end{equation}
\end{definition}
%%%%%%%%%%%%%%%%%%%%%

This definition generalizes to multivariate functions by requiring GSC along every straight line. Specifically, let $h:\Rn\to\R$ be a closed convex, lower semi-continuous function with open and convex effective domain $\dom h=\scrC\subset\Rn$. For $x\in \scrC$ and $u,v\in\Rn$, define the real-valued function $\phi(t):=\inner{\nabla^{2}h(x+tv)u,u}$. For $t\in\dom\phi$, one sees that 
\begin{align*}
\phi'(t)=\inner{D^{3}h(x+tv)[v]u,u}, 
\end{align*}
so that we can define generalized self-concordance of a function by formulating conditions on the behavior of $\phi'(0)$.
%%%%%%%%%%%%%%%%%%%%%%%%%%%%%%%
\begin{definition}\citep{SunTran18}
A closed convex function $h\in\bC^{3}(\dom h)$, with $\dom h$ open, is called $(M,\nu)$ generalized self-concordant of the order $\nu>0$ and constant $M\geq 0$ if for all $x\in \dom h$ 
\begin{equation}
\abs{\inner{D^{3}h(x)[v]u,u}}\leq M\norm{u}^{2}_{x}\norm{v}^{\nu-2}_{x}\norm{v}^{3-\nu}_{2}\qquad\forall u,v\in\Rn. 
\end{equation}
We denote this class of functions as $\scrH_{M,\nu}(\dom h)$. 
\end{definition}
%%%%%%%%%%%%%%%%%%%%%%%%%%%%%%%

As in the theory of standard self-concordant functions, the precise value of the scale parameter $M>0$ is not of big importance for theoretical considerations. In fact, it is easy to see that we can always rescale the function so that the definition of a GSC function holds for $M=2$.
\begin{lemma}
If $\phi\in\bC^{3}(\dom\phi;\R)$ belongs to the class $\scrH_{M,\nu}(\dom h)$, then $\left(\frac{M}{2}\right)^{\frac{2}{\nu-2}}\phi\in\scrH_{2,\nu}(\dom h)$.
\end{lemma}
\begin{Proof}
Let $\psi(t):=\left(\frac{M}{2}\right)^{\frac{2}{\nu-2}}\phi(t)$. Then, for all $t\in\dom\phi$, we have 
\begin{align*}
\abs{\psi'''(t)}&=\left(\frac{M}{2}\right)^{\frac{2}{\nu-2}}\abs{\phi'''(t)}\leq M\left(\frac{M}{2}\right)^{\frac{2}{\nu-2}}\phi''(t)^{\nu/2}=2\psi''(t)^{\nu/2}.
\end{align*}
\end{Proof}
The function $h\in\scrH_{M,\nu}(\dom h)$ defines a semi-norm 
\begin{equation}\label{eq:RM}
\norm{d}_{x}:=\sqrt{\inner{\nabla^{2}h(x)d,d}},
\end{equation}
with dual norm 
\begin{equation}
\norm{d}^{\ast}_{x}:=\sup_{d\in\Rn}\{2\inner{d,a}-\norm{d}_{x}^{2}\}.
\end{equation}
Note that if $H(x)\equiv\nabla^{2}h(x)\succ 0$ then $\norm{\cdot}_{x}$ is a real norm (necessarily equivalent to the euclidean norm), and $\norm{d}^{\ast}_{x}=\sqrt{\inner{[H(x)]^{-1}d,d}}$. The barrier-character of functions $h\in\scrH_{M,\nu}(\dom h)$ is made clear in the next Lemma. 

\begin{lemma}[\cite{NesConvex},Thm. 4.1.4]
\label{lem:barrier}
For every sequence $(x^{k})_{k\geq 0}$ such that $(x^{k})_{k\geq 0}\subset\dom h$ and $x^{k}\to x\in\bd(\dom h)$, we have $\lim_{k\to\infty}h(x^{k})=\infty$.
\end{lemma}
\begin{Proof}
For all $k\geq 0$, we have 
\begin{align*}
h(x^{k})\geq h(x^{0})+\inner{\nabla h(x^{0}),x^{k}-x^{0}}. 
\end{align*}
If the sequence $(h(x^{k}))_{k\geq 0}$ is bounded from above, we can descent to a subsequence along which $h(x^{k})\to \bar{h}$ (we omit the relabeling). Then, for all $k\geq 0$, 
$z^{k}=(x^{k},h(x^{k}))\in \epi(h)$, and $z^{k}\to z=(x,\bar{h})\in\epi(h)$, since the function is closed. Hence, $x\in \dom h$. A contradiction. 
\end{Proof}

Given $\nu\in(2,4]$ and $h\in\scrH_{M,\nu}(\scrC)$, we define the distance function
\begin{equation}
\label{eq:distance}
\metric_{\nu}(x,y):=\left\{\begin{array}{lr} 
M\norm{y-x}_{2} & \text{ if }\nu=2,\\
\frac{\nu-2}{2}M\norm{y-x}^{3-\nu}_{2}\cdot\norm{y-x}^{\nu-2}_{x} & \text{if }\nu>2. 
\end{array}\right.
\end{equation}
The \emph{Dikin Ellipsoid} with respect to the distance function $\metric_{\nu}$ is defined as 
\begin{equation}\label{eq:dikin}
\scrW(x;r):=\{y\in\Rn: \metric_{\nu}(x,y)<r\}\quad\forall (x,r)\in\dom h\times\R.
\end{equation}
%The notion of generalized self-concordance allows us to induce feasible numerical schemes by controlling the step-length of the method ensuring that the subsequent iterates are confined to stay within a Dikin ellipsoid of radius smaller than 1. This is sufficient to keep the process in the interior of the feasible set, thanks to the next lemma. 

\begin{lemma}[\cite{SunTran18}, Prop. 7]
\label{lem:Dikin}
Let $h\in\scrH_{M,\nu}(C)$ be a barrier-generating kernel of order $\nu\in(2,4]$. We have $\scrW(x;1)\subseteq\dom h $ for all $x\in\dom h$.
\end{lemma} 

\begin{remark}
The familiar inclusion $\scrW(x;1)\subset\dom h$ is only true if $\nu>2$. This is intuitive, since for $\nu=2$, the local norm effectively boils down to the euclidean norm, and thus is not adaptive to the local geometry. As a consequence, our algorithmic scheme will take as inputs functions $h\in\scrH_{M,\nu}(\scrC)$ with $\nu>2$ only. This covers the important case of standard self-concordant functions, as well as many entropy-based barrier functions familiar from the literature on Bregman proximal gradient methods. However, our method also works well for generalized self-concordant function of order $\nu\in(3,4]$. This range cannot be analyzed by proximal based algorithms studied in recent work \cite{CevKyrTra15,SunTran18} and \cite{LuSunTran19}.
\end{remark}

We define the \emph{Bregman divergence} associated to $h\in\scrH_{M,\nu}(\scrC)$ as
\begin{equation}\label{eq:D}
D_{h}(x,y):=h(x)-h(y)-\inner{\nabla h(y),x-y}\quad\text{for } x\in\scrC,y\in\scrC.
\end{equation}
Since this divergence function will be a crucial quantity of interest in measuring the per-iteration progress of our method, it is instrumental to have universal bunds on the function values. For the class of self-concordant functions, such bounds are classical to the field (see e.g. \cite{NesNem94}). For the  Bregman divergence induced by the class of generalized self-concordant functions, a similar universal bound can be reported. 
\begin{lemma}[\cite{SunTran18}, Prop. 10]
\label{lem:Bregman}
Let $x\in\dom h$ for $h\in\scrH_{M,\nu}(\scrC)$ and $\nu\in(2,4]$. Then 
\begin{equation}\label{eq:Bregbound}
\omega_{\nu}(-\metric_{\nu}(x,y))\norm{y-x}^{2}_{x}\leq D_{h}(y,x)\leq \omega_{\nu}(\metric_{\nu}(x,y))\norm{y-x}^{2}_{x},
\end{equation}
for all $y\in\scrW(x;1)$, where  
\begin{equation}
\omega_{\nu}(t):=\left\{\begin{array}{ll}
%\frac{1}{t^{2}}(e^{t}-t-1) & \text{if }\nu=2,\\
\frac{-t-\ln(1-t)}{t^{2}} & \text{if }\nu=3,\\
\frac{(1-t)\ln(1-t)+t}{t^{2}} & \text{ if }\nu=4,\\
\left(\frac{\nu-2}{4-\nu}\right)\frac{1}{t}\left[\frac{\nu-2}{2(3-\nu)t}((1-t)^{\frac{2(3-\nu)}{2-\nu}}-1)-1\right] & \text{ otherwise.}
\end{array}\right. 
\end{equation}
\end{lemma}
%%%%%%%%%%%%
Let $h\in\scrH_{M,\nu}(\scrC)$ and $x_{j}\to \bd(\scrC)=\bar{\scrC}\setminus \scrC$. Then, by Lemma \ref{lem:barrier}, $h(x_{j})\to\infty$. We claim that this implies $\norm{\nabla h(x_{j})}_{2}\to\infty$. Indeed, by convexity, for all $y\in \scrC$, we have 
\begin{align*}
h(y)\geq h(x_{j})+\inner{\nabla h(x_{j}),y-x_{j}}. 
\end{align*}
Therefore, by Cauchy-Schwarz  
\begin{align*}
h(y)\geq h(x_{j})-\norm{y-x_{j}}\cdot\norm{\nabla h(x_{j})}_{2}
\end{align*}
If $\norm{\nabla h(x_{j})}_{2}$ would be bounded, the right-hand side diverges to $\infty$, whereas the left hand side is bounded for $y\in\dom (h)$. This gives a contradiction. We conclude that $h$ is \emph{Legendre:} 
\begin{equation}\label{eq:legendre}
x_{j}\to x^{\ast}\in \bd(\scrC)\Rightarrow \norm{\nabla h(x_{j})}_{2}\to\infty. 
\end{equation} 
%Based on this, we can deduce the following result which coincides with \cite[Lem. 4.2]{ABB04}. 

%\begin{lemma}\label{lem:normal}
%Let $(x_{j})_{j\geq 1}\subset \scrC$ be such that $x_{j}\to x^{\ast}\in\bd(\scrC)=\bar{\scrC}\setminus \scrC$. and $\frac{\nabla h(x_{j})}{\norm{\nabla h(x_{j})}}\to p\in\Rn$ with $h\in\scrH_{M,\nu}(\scrC)$. Then $p\in \NC_{\bar{\scrC}}(x^{\ast})$.
%\end{lemma}
%\begin{proof}
%By convexity of $h$, $\inner{\nabla h(x_{j})-\nabla h(y),x_{j}-y}\geq 0$ for all $y\in \scrC$. Dividing by $\norm{\nabla h(x_{j})}$ and letting $j\to\infty$, we get $\inner{p,y-x^{\ast}}\leq %0$ for all $y\in \scrC$, which holds also for $y\in\bar{\scrC}$. Hence, $p\in \NC_{\bar{\scrC}}(x^{\ast})$. 
%\end{proof}

\begin{lemma}\label{lem:nondeg}
If $h\in\scrH_{M,\nu}(\dom h)$ and $\dom h\subset\Rn$ contains no lines, then $H(x)\succ 0$ for all $x\in\dom h$. 
\end{lemma}
\begin{Proof}
Define the recessive subspace $E_{h}:=\{d\in\Rn:\norm{d}_{x}=0\text{ for some }x\in\dom h\}$. From \cite[Prop. 8]{SunTran18} we deduce that for all $r\in(0,1)$ and $y\in\Rn$ with $r=\metric_{\nu}(x,y)$, we have 
\begin{align}\label{eq:Hessianbound}
(1-r)^{\frac{2}{\nu-2}}\nabla^{2}h(x)\leq \nabla^{2}h(y)\leq(1-r)^{\frac{-2}{\nu-2}}\nabla^{2}h(x). 
\end{align}
Let $Z_{d}:=\{x\in\dom h: \norm{d}_{x}=0\}$. \eqref{eq:Hessianbound} implies that $x\in Z_{d}\Rightarrow y\in Z_{d}$ for all $y\in \scrW(x;r)$, and therefore $Z_{d}$ is open. Since $h\in\bC^{3}(\dom h)$, it is closed as well. Therefore $Z_{d}$ is either empty or the entire set $\Rn$. This implies that $E_{h}$ is either empty or $\Rn$. From here the result follows from \cite[Thm. 4.1.3]{NesConvex}.
\end{Proof}

In order to fully understand the behavior of Newton methods involving standard self-concordant function $h\in\scrH_{2,3}(\dom h)$, the general theory laid out by \cite{NesNem94} introduced the concept of a \emph{self-concordant barrier} (SC-B). 
\begin{definition}\label{def:SCB}
For some scalar $\theta\geq 0$, a function $h\in\scrH_{2,3}(\scrC)$ is a self-concordant barrier of order $\theta>0$ ($\theta$-SCB) if 
\begin{equation}\label{eq:SC_barrier}
\norm{\nabla h(x)}_{x}^{\ast}\leq\sqrt{\theta}\qquad\forall x\in\scrC.
\end{equation}
\end{definition}
It is very remarkable that every convex body, i.e. every compact convex set with nonempty interior, admits a $\theta$-SCB with $\theta$ a dimension-dependent constant. This class of functions is the main driver in standard interior-point solvers, and conceptually it reveals the two main problems IPMs face when confronted with large-scale problems: (i) Even if a suitable self-concordant barrier can be computed, the order parameter $\theta$ is dimension dependent, and thus iteration complexity cannot be dimension-free. However, Bubeck and Eldan \cite{BubEld18} have recently shown that the Nesterov-Nemirovski universal barrier (see below) is a $\theta$-SCB on $\scrC$ with $\theta=(1+\eps_{n})n$ and $\eps_{n}\leq 100\sqrt{\log(n)/n}$.  (ii) The construction of the \emph{universal barrier} for a convex body $\bar{\scrC}\subset\Rn$ due to Nesterov and Nemirovski \cite{NesNem94}  is based on the log-Laplace transform $h^{\ast}(w)=\log\left(\int_{\bar{\scrC}}\exp(\inner{x,w})\dif x\right)$ and its convex conjugate $h(x)=\sup_{w\in\Rn}\{\inner{w,x}-h^{\ast}(w)\}$. Unless the set $\scrC$ is special, computing a universal barrier is infeasible. 

%From the point of view of the theory of first-order methods, our implicit assumption is thus similar to the prevailing implicit assumption behind Bregman proximal algorithms. In the latter family of algorithms, it is implicitly assumed that the proximity operator can be evaluated efficiently, which is also essentially an assumption on the convex geometry of the feasible set. As will be seen, the advantage of the Hessian-barrier approach is that we do not have to evaluate a proximity operator at each iteration of the algorithm. Instead, the search directions \acl{HBA} uses are always available in closed-form expressions and their computations involve the solution of a linear system of equations. We believe that this is a significant advantage compared to Bregman proximal algorithms. 

%%%%%%%%%%%%%%%%%%%%%%%%%%%%%%%%%%%%%%%%%%%
%----------------------------------------------------------------------
%%% Kernels
%----------------------------------------------------------------------
\subsection{Barrier generating kernels}
\label{sec:bgk}

The class of metric generating functions of interest in this paper is defined as a subset $\scrF_{M,\nu}(\scrC)$ contained in $\scrG(\scrC)\cap\scrH_{M,\nu}(\scrC)$, whose exact definition depends on the order parameter $\nu$. If $\nu\in(2,3]$, we take $\scrF_{M,\nu}(\scrC)=\scrG(\scrC)\cap\scrH_{M,\nu}(\scrC)$. If $\nu\in (3,4]$, we additionally assume \emph{coercivity} of the kernel-generating distance $h$, i.e that $h(x)\to\infty \text{ whenever }\norm{x}\to\infty.$
\begin{definition}\label{def:bgk}
The class of \emph{barrier-generating kernels} is defined as
\begin{align*}
\scrF_{M,\nu}(\scrC):=\left\{\begin{array}{lr} 
\scrG_{M,\nu}(\scrC)\cap\scrH_{M,\nu}(\scrC)\text{ if }\nu\in(2,3],\\
\{h\in\scrG_{M,\nu}(\scrC)\cap\scrH_{M,\nu}(\scrC): h\text{ is coercive}\} & \text{if }\nu\in(3,4].
\end{array}\right.
\end{align*}
\end{definition}

\subsection{Examples}
\label{ex:bgf}
In order to illustrate the flexibility of the framework of sets endowed with barrier-generating kernels, we collect below some representative examples taken from the literature. For many more examples, we refer the reader to \cite{NesNem94}. 

The first set of examples are tailored to product domains $\scrC=\scrC_{1}\times\cdots\times\scrC_{n}$, where each $\scrC_{i}$ is an open convex subset of the real-line. For such domains, \emph{decomposable} barrier-generating kernels are an attractive choice since their Hessian matrix is diagonal. Specifically, we consider functions $h\in\scrF_{M,\nu}(\scrC)$ of the form 
\[
h(x_{1},\ldots,x_{n})=\sum_{i=1}^{n}\phi_{i}(x_{i}), 
\]
where each function $\phi_{i}\in\scrF_{M_{\phi_{i}},\nu}(\scrC_{i})$. By \cite[Prop. 1]{SunTran18}, $h$ is generalized self-concordant with $\dom h =\bigcap_{i=1}^{n}\dom\phi_{i}$ and constant $M:=\max\{M_{\phi_{1}},\ldots,M_{\phi_{n}}\}$. For different structures of the sets $\scrC_{i}$, we can propose different barrier-generating kernels $\phi_{i}$. Here are some illustrative examples.

%The next examples show that this intersection contains a plethora of important kernel generating distances, well known in the theory of first-order methods. Each example is a one-dimensional $\varphi$, belonging to the class $\scrF_{M,\nu}(\dom\varphi)$ for some $\nu>2$. To obtain a corresponding function $h\in\scrF_{M_{h},\nu}(\dom h)$, simply use the formula $h(x)=\sum_{i=1}^{n}\varphi_{i}(x_{i})$, where $\varphi_{i}\in\scrF_{M_{\varphi_{i}},\nu}(\dom\varphi_{i})$. 

% In order the verify the continuity of the dual norm $\norm{\nabla h(x)}^{\ast}_{x}=\left(\sum_{i=1}^{n}\frac{\varphi'_{i}(x_{i})^{2}}{\varphi''_{i}(x_{i})}\right)^{1/2}$ (required by Definition \ref{def:kgd}(c)), it suffices to verify that the functions $t\mapsto \frac{\varphi'_{i}(t)^{2}}{\varphi''_{i}(t)}$ are continuous. 

\begin{enumerate}
\item \emph{Burg entropy:} $\phi(t)=-\log(t)$ for $t>0$ is an element of $\scrF_{2,3}(\R_{++})$; %Note that $\frac{\varphi'(t)^{2}}{\varphi''(t)}=1$ for all $t>0$.
\item \emph{Entropy-Barrier:} $\phi(t)=t\log(t)-\log(t)$ for $t>0$ is an element of $\scrF_{2,3}(\R_{++})$;
\item Consider the function $\phi(t)=(1-t/\kappa)^{-\kappa}$ for $\kappa>0$ and $t\in(-\infty,\kappa)$. Then $\dom\phi=(-\infty,\kappa)$ and one can check that $\phi\in\scrF_{M,\nu}((-\infty,\kappa))$, where $M=\frac{2+\kappa}{\kappa}\left(\frac{\kappa}{1+\kappa}\right)^{\frac{1}{2+\kappa}}$ and $\nu=\frac{2(3+\kappa)}{2+\kappa}\in(2,3)$ for $\kappa>0$.  
\item The function $\phi(t)=\frac{1}{\sqrt{1-t^{2}}}$ defines an element in $\scrF_{M,\nu}((-1,1))$ for $M$ a constant smaller than 3.25 and $\nu=14/5$. 
\end{enumerate}

\begin{remark}\label{rem:box}
Let $\bar{\scrC}=\prod_{i=1}^{n}[a_{i},b_{i}]$ be a high-dimensional box of dimension $n\gg 1$ where $+\infty>b_{i}\geq a_{i}>-\infty$. According to \cite{CoxJudNem14}, this geometry is computational challenging for standard proximal methods.\footnote{To be clear, the challenge is not to find a good distance generating function, but rather the scalability of the mirror descent algorithm. We refer the reader to in-depth discussion in \cite{CoxJudNem14} for details.} This geometry can be easily endowed with a barrier-generating kernel $h\in\scrF_{M,\nu}(\scrC)$ given by a sum of Burg entropies or Entropy barriers, for instance, leading to a simple Riemannian metric on the interior of this box.
\end{remark}

The above examples provide a snapshot of common barrier-generating kernels used in practice. However, it is possible to combine these functions to obtain mixture functions that preserve the properties imposed on an element $h\in\scrF_{M,\nu}(\scrC)$. In particular, it is easy to find barrier generating kernels for geometries which are given as intersections of open convex sets $\scrC_{1},\ldots,\scrC_{J}$, each admitting a generalized self-concordant function $h_{j}\in\scrF_{M_{j},\nu}(\scrC_{j}),1\leq j\leq J$. This appears in second-order cone programming problems, which have as a special case optimization problems with quadratic constraints (see \cite{LobVanBoy98} for a survey). The typical sets appearing in such optimization problems are the following: 

\begin{itemize}
\item $\scrC=\Rn_{++},h(x)=-\sum_{i}\log(x_{i})$ is a $n$-SCB and an element of $h\in\scrF_{2,3}(\scrC)$;
\item Consider the second-order cone $\bar{\scrC}\equiv \scrL^{n}:=\cl \left(\{x=(t,w)\in\R\times \R^{n-1}: t> \norm{w}_{2}\}\right)$.
The function 
\begin{align*}
h(x)=-\log(t^{2}-\norm{w}_{2}^{2})\quad x:=(t,w)\in \scrC,
\end{align*}
is a barrier-generating kernel belonging to $\scrF_{2,3}(\scrC)$. It is also an $2$-SCB. 
\item Consider the cone of positive definite symmetric $n\times n$ matrices with real entries $\Sigma^{n}_{++}:=\{x\in\R^{n\times n}:x\succ 0,x^{\top}=x\}$, and set
$\bar{\scrC}:=\cl(\Sigma^{n}_{++})$. The function $h(x)=-\log\det(x)$ for $x\in \scrC$ is a barrier-generating kernel of class $\scrF_{2,3}(\scrC)$. It is also an $n$-SCB. 
\item Let $B$ be a $p\times n$ matrix with rows $b^{\top}_{1},\ldots,b^{\top}_{p}$, and $d$ a given vector in $\R^{p}$. Consider the polyhedral set $\bar{\scrC}:=\{x\in\Rn: Bx\leq d\}.$
Assume that $\scrC=\{x\in\Rn: Bx<d\}$ is nonempty (Slater condition). Then, the function 
$h(x)=\sum_{j=1}^{p}-\log(d_{j}-b^{\top}_{j}x)$ is a barrier generating kernel belonging to the class $\scrF_{2,3}(\scrC)$.
\item Let $\scrM_{m,n}$ be the space of real $m\times n$ matrices with inner product $\inner{A,B}=\tr(AB^{\top})$. The standard operator norm is defined as $\abs{Q}:=\max\{\norm{Qw}:\norm{w}=1\}$. Consider the set $\bar{\scrC}=\cl\left(\{x=(t,Q)\in\R\times\scrM_{m,n}: t>\abs{Q}\}\right).$
This set admits a barrier-generating kernel $h\in\scrF_{2,3}(\scrC)$ given by $h(x)=\log\det(t\Id-\frac{1}{t}QQ^{\top})-\log(t)$ for $x=(t,Q)\in\scrC.$
\end{itemize}
%------------------------------
%%%%% Optimization Problem %%%%%
%------------------------------
\subsection{The minimization problem}
We are given a matrix $A\in\R^{m\times n}$ of full row rank $m$ and $b\in\image(A)$. Define the sets $\scrA=\{x\in\Rn: Ax=b\}$, and $\scrA_{0}=\ker(A)$, so that $\scrA_{0}^{\bot}=\image(A^{\top})$. Let $\scrC$ be a nonempty open convex set in $\Rn$ with closure $\bar{\scrC}$ that is not contained in any $(n-1)$-dimensional affine subspace. Throughout the rest of this paper the following assumption is taken as a standing hypothesis. 
\begin{assumption}
\label{ass:C}
The set $\scrC$ is nonempty, convex and contains no lines. 
\end{assumption}
Combining this assumption with Lemma \ref{lem:nondeg}, we know that the Hessian matrix $H(x)=\nabla^{2}h(x)$ is positive definite on $\scrC$. The matrix-valued function $H:\scrC\to\Sigma^{n}_{++}$ defines a Riemannian manifold $(\scrC,\norm{\cdot}_{x})$, with Riemannian metric given by \eqref{eq:RM}. We are also given a lower semi-continuous $f:\Rn\to(-\infty,+\infty]$. The problem we aim to solve is the minimization problem 
\begin{equation}\tag{P}
f^{\ast}:=\inf\{f(x): x\in\bar{\scrC},Ax=b\}
\end{equation}
The feasible set of \eqref{eq:P} is denoted as $\scrX=\bar{\scrC}\cap\scrA$, and we shall denote by $\scrX^{\circ}$ the relative interior of $\scrX$, that is, $\scrX^{\circ}=\{x\in\Rn: x\in \scrC,Ax=b\}.$ As a standing hypothesis, we shall impose the following Slater constraint qualification condition: 
%%%%%%%%%%%%%%%%%%%%%%%%
\begin{assumption}\label{ass:slater}
$\scrX^{\circ}\neq\emptyset$.
\end{assumption}
%%%%%%%%%%%%%%%%%%%%%%
%\begin{assumption}\label{ass:fbounded}
%The objective function is bounded from below in the feasible set.
%\end{assumption} 
%%%%%%%%%%%%%%%%%%%%%
\begin{assumption}\label{ass:levelbound}
The level sets of the objective function are bounded: Given $x^{0}\in\scrX^{\circ}$ there exists $R>0$ such that $\sup\{\norm{x}_{\infty}: f(x)\leq f(x^{0})\}\leq R$.
\end{assumption} 
%%%%%%%%
For $\eps>0$, an $\eps$ global minimizer is defined as a feasible solution $x_{\eps}$ such that 
\begin{equation}
f(x_{\eps})-\inf_{x\in\scrX}f(x)\leq\eps. 
\end{equation}
It is well known that finding an $\eps$-global minimizer is a strongly NP-hard problem (see e.g. \cite{GeJiaYe11}). Even worse, it is also well known that in general, finding a descent direction for a non-convex non-smooth function is NP-hard. As concrete illustration, even deciding whether the function 
\begin{align*}
f(x)=(1-1/\gamma)\max_{1\leq i\leq n}\abs{x_{i}}-\min_{1\leq i\leq n}\abs{x_{i}}+\abs{\inner{c,x}},
\end{align*}
where $x\in\Rn,c\in\N^{n}$ and $\gamma=\sum_{i=1}^{n}c_{i}$, admits a descent direction is NP-hard \cite[Lem. 1]{Nes13}. Therefore, in this paper we restrict ourselves to objective functions of very special structure. Namely, we consider the problem of minimizing a real-valued function $f:\Rn\to(-\infty,\infty]$ which is continuously differentiable on an open convex set $\scrC\subset\Rn$, and possibly non-differentiable at the boundary $\bd\scrC=\bar{\scrC}\setminus\scrC$. 
%%%%%%%%%%%%%%%%%%%%%%%%%%%%%%%%%
\begin{assumption}\label{ass:fcontinuous}
$f:\Rn\to(-\infty,\infty]$ is a proper and lower semi-continuous function with $f\in\bC^{1}(\scrC)$.
\end{assumption} 
%%%%%%%%%%%%%%%%%%%%%%%%
%The first-order optimality conditions at the point of local minimum $x^{\ast}\in\scrX$ at which $f$ is directionally differentiable can be written as 
%\begin{equation}
%f'(x^{\ast};d)\geq 0\qquad\forall d\in\TC_{\scrX}(x^{\ast})=\TC_{\bar{C}}(x^{\ast})\cap\scrA_{0},
%\end{equation}
%where $f'(x;d)$ denotes the directional derivative of $f$ at $x^{\ast}$ in direction $d$. This last variant of first-order optimality is convenient for defining an approximate solution:
%\begin{definition}
%\PD{Do we need this definition? We do not use it as far as I see.} The point $\bar{x}\in\scrX^{\circ}$ is an $\eps$-optimal solution if 
%\begin{equation}
%f'(\bar{x};d)\geq -\eps\qquad\forall d\in\TC_{\scrX}(x^{\ast}). 
%\end{equation}
%\end{definition}
%Note that for $\bar{x}\in\scrX^{\circ},$ the set of feasible directions reduces to  
%\begin{equation}
%\inner{\nabla f(\bar{x}),d}\geq-\eps\qquad\forall d\in\scrA_{0}.
%\end{equation}

The smoothness condition formulated in Assumption \ref{ass:fcontinuous} is silent about the behavior of the function at the boundary $\bd(\scrC)$. In case where the function $f$ is twice continuously differentiable Cartis, Gould and Toint defined in \cite{CarGouToi12} the following criticality measure at $x\in\scrX$ 
\begin{equation}
\chi_{CGT}(x):=\abs{\min_{x+d\in\scrX,\norm{d}_{2}\leq 1}\inner{\nabla f(x),d}}.
\end{equation}
They subsequently proved $\bigoh(\eps^{-2})$ iteration complexity for reaching a point with $\chi_{CGT}(x)\leq\eps$. We propose a similar criticality measure here, but make use of the local norm. In particular, we consider the primal-dual stationarity measure at $(x,y)\in\scrX^{\circ}\times\R^{m}$ given by 
\begin{equation}\label{eq:chi}
\chi(x,y):=\norm{\nabla f(x)-A^{\top}y}^{\ast}_{x}.
\end{equation}
%%%%%%%%%%%%%%%%
\begin{definition}\label{def:stationary}
Given $\eps>0$, a pair $(x^{\ast},y^{\ast})\in\scrX^{\circ}\times\R^{m}$ is calld $\eps$-stationary if $\chi(x,y)\leq\eps$.
\end{definition}
%%%%%%%%%%%%%%%%%%

In order to motivate this criticality measure, we first recall the classical Fenchel-Young inequality 
\begin{equation}
\abs{\inner{u,v}}\leq\norm{u}_{x}\cdot\norm{v}^{\ast}_{x}\qquad\forall u,v\in\Rn,x\in\scrX^{\circ}. 
\end{equation}
Hence, for $v=\nabla f(x)-A^{\top}y$ and $u\in\scrA_{0}$ with $\norm{u}_{x}=1$, this inequality readily gives us 
\[
-\chi(x,y)\leq\inner{\nabla f(x)-A^{\top}y,u}=\inner{\nabla f(x),u}\leq\chi(x,y),
\]
and in particular,
\[
\abs{\min_{u\in\scrA_{0},\norm{u}_{x}=1}\inner{\nabla f(x),u}}\leq\chi(x,y). 
\]
Thus, the primal-dual criticality measure is an upper bound of a version of the Cartis-Gould-Toint criticality measure $\chi_{CGT}(x)$, and we note in passing that if $\chi(x,y)\leq \eps$, then automatically $\inner{\nabla f(x),u}\in(-\eps,\eps)$ for all $u\in\scrA_{0}$ satisfying $\norm{u}_{x}=1$. One potentially troublesome part in the definition of our proposed criticality measure is that it is formulated in terms of the local norm. Hence, we would need to evaluate the inverse matrix $[H(x)]^{-1}$ (provided it exists at $x$). However, for our algorithm this is not a problem since we will have very good control about the location of the iterates. Indeed, as will be seen in Section \ref{sec:complexity}, the algorithm will take values on a compact set $\scrS_{\mu}(x^{0})$ in $\scrX^{\circ}$ (this is similar to proximal based self-concordant algorithms and exploited in \cite{Lu17} in the convergence analysis). On this set, we have very good control on the eigenvalues of the Hessian matrix $H(x)=\nabla^{2}h(x)$, and in fact, under assumption spelled out explicitly in the sections to follow, we can provide upper and lower bounds on the eigenvalues of $H(x)$ over the set $\scrS_{\mu}(x^{0})$, denoted as $0<\sigma_{h}<\tau_{h}<\infty$. Therefore, during the working phase of the algorithm, we produce a primal-dual sequence $(x^{k},y^{k})$ along which the criticality measure is sandwiched as 
\begin{equation}\label{eq:sandwich1}
\tau_{h}^{-1/2}\norm{\nabla f(x^{k})-A^{\top}y^{k}}_{2}\leq\norm{\nabla f(x^{k})-A^{\top}y^{k}}^{\ast}_{x^{k}}\leq\sigma_{h}^{-1/2}\norm{\nabla f(x^{k})-A^{\top}y^{k}}_{2}. 
\end{equation}
Therefore, if the euclidean norm of the vector $\nabla f(x^{k})-A^{\top}y^{k}$ falls below a cut-off $\eps>0$, we have reached an $\eps$-stationary point in the sense of Definition \ref{def:stationary}.

\begin{remark}
This notion of stationarity is also motivated by the structure of the KKT conditions satisfied by a solution candidate for problem \eqref{eq:P} taking values in the relative interior $\scrX^{\circ}=\scrC\cap\scrA$. As a concrete illustration, let us consider the set $\bar{\scrC}=\Rn_{+}$, so that we are in the setting of \cite{HBA-linear}. The complementary slackness condition for the resulting optimization problem \eqref{eq:P} reads as 
\begin{equation}
X(\nabla f(x)-A^{\top}y)=0
\end{equation}
where $X=\diag\{x_{1},\ldots,x_{n}\}$. Hence, a reasonable definition of an $\eps$-KKT point under a Riemannian-Hessian structure induced by the Hessian of the function $h(x) =- \sum_{i=1}^n \ln x_i$ would read as 
\begin{align*}
\norm{X(\nabla f(x)-A^{\top}y)}_{\infty}\leq \norm{\nabla f(x)-A^{\top}y}^{\ast}_{x}\leq \eps. 
\end{align*}
This $\eps$-KKT definition has also been used in \cite{HaeLiuYe18}.
\end{remark}

%% file: HBA_general.tex
%%% HBA_general
%----------------------------------------------------------------------
% !TEX root = ./Main.tex
%
In this section we describe a conceptual version of the Hessian-barrier method. To this end, we are given an open nonempty set $\scrC\subset\Rn$ satisfying Assumption \ref{ass:C}, admitting a computable barrier generating kernel $h\in\scrF_{M,\nu}(\scrC)$.

\subsection{Defining the search directions}
For a pair $(x,g)\in \scrC\times\Rn$, define the functions
\begin{align*}
\psi(x,g)&:=\min_{v}\left\{\inner{g,v}+\frac{1}{2}\norm{v}_{x}^{2}: Av=0\right\},\text{ and }\\
V(x,g)&:=\argmin_{v}\left\{\inner{g,v}+\frac{1}{2}\norm{v}_{x}^{2}: Av=0\right\}.
\end{align*}
Computing the vector $V(x,g)$ means finding a pair $(V(x,g),y(x,g))=(v,y)\in\Rn\times \R^{m}$ solving the Newton-type of system
\begin{equation}\label{eq:finder}
\left[\begin{array}{cc} 
H(x) & -A^{\top}\\
-A & 0
\end{array}\right]\cdot\left[\begin{array}{c} v \\ y\end{array}\right]=\left[\begin{array}{c} -g \\ 0 \end{array}\right]. 
\end{equation}
In particular, the complexity of computing $V(x,g)$ is of the same order as finding a Newton direction, and the practical efficiency of the method depends heavily on the structure of the matrices $H(x)$ and $A$, respectively. In any case, given that $H(x)\in\Sigma^{n}_{++}$, we obtain a closed form expression for the vector $V(x,g)$ as 
\begin{align}\label{eq:V}
V(x,g)&=-P_{x}[H(x)]^{-1}g\qquad\forall (x,g)\in \scrC\times\Rn,\text{ and }\\
y(x,g)&=(A[H(x)]^{-1}A^{\top})^{-1}A[H(x)]^{-1}g
\label{eq:dual}
\end{align}
where the matrix valued function $P:\scrC\to \R^{n\times n}$ defined in \eqref{eq:project}.
%\begin{equation}
%P_{x}=\Id-[H(x)]^{-1}A^{\top}(A[H(x)]^{-1}A^{\top})^{-1}A.
%\end{equation}
We just remark that, given the matrix $A$ being of full rank, the function $x\mapsto A[H(x)]^{-1}A^{\top}$ is invertible \cite{HBA-linear}. Computational efficiency considerations will be made later. We close this section by establishing some general properties of the mapping $\psi$ and $V$.

\begin{proposition}\label{prop:propertiesV}
The following assertions are true:
\begin{itemize}
\item[(a)] The mapping $\psi(x,\cdot):\Rn\to\R$ is concave and continuously differentiable with 
\[
\nabla_{g}\psi(x,g)=V(x,g)\qquad\forall (x,u)\in \scrC\times\Rn.
\]
\item[(b)] If $h$ is $K$-strongly convex under the $\ell_{2}$ norm, then for every $x\in \scrC$, the mappings $V(x,\cdot)$ and $\psi(x,\cdot)$ are $\frac{1}{K}$-Lipschitz,

\item[(c)] If $g=p+z\in\NC_{\scrX}(x)$ with $p\in\NC_{\bar{\scrC}}(x),z\in\scrA_{0}^{\bot}$, then 
\[
V(x,p+z)=V(x,p)
\]
\end{itemize}
\end{proposition}
\begin{Proof}
Since $\psi(x,g)$ is the pointwise minimum of a linear function, it must be concave. The integrability condition on the vector field $g\mapsto V(x,g)$ is a straightforward computation. Parts (b) and (c) are standard, and follow from the general analysis of such projection schemes as in \cite{Nes05}. It is however instructive here to go over the computations. First, the $K$-strong convexity of the norm ensures that $u\mapsto \psi(x,u)$ is well defined and convex. In particular, $g\mapsto V(x,g)$ is uniquely defined by eq. \eqref{eq:V}. Therefore,  
\[
\psi(x,g)=\frac{1}{2}\norm{V(x,g)}^{2}_{x}\qquad\forall (x,g)\in \scrC\times\Rn.
\]
Let $g_{1},g_{2}\in\Rn$ be arbitrary and set $v_{1}=V(x,g_{1}),v_{2}=V(x,g_{2})$. The optimality conditions at a given point $x\in \scrC$ imply
\begin{align*}
&(g_{1}-H(x)v_{1})^{\top}(v_{2}-v_{1})=0 \text{ and}\\
&(g_{2}-H(x)v_{2})^{\top}(v_{1}-v_{2})=0.
\end{align*}
Adding both, and using the $K$-strong convexity shows 
\begin{align*}
(g_{1}-g_{2})^{\top}(v_{1}-v_{2})&=(H(x)(v_{1}-v_{2}))^{\top}(v_{1}-v_{2})\geq K\norm{v_{1}-v_{2}}^{2}_{2}.
\end{align*}
Using the Cauchy-Schwarz inequality, we arrive at 
\begin{align*}
\norm{V(x,g_{1})-V(x,g_{2})}_{2}\leq\frac{1}{K}\norm{g_{1}-g_{2}}.
\end{align*}
\end{Proof}
%%%%%%%%%%%%%%%%%%%%%%%%%%%%%
\begin{corollary}
For all $x\in\scrX^{\circ}$, we have 
\begin{align*}
V(x,g)=0\qquad\forall  g \in \NC_{\scrX}(x)=\NC_{\bar{\scrC}}(x)+\scrA^{\bot}_{0}.
\end{align*}
\end{corollary}
\begin{Proof}
Just observe that for $x\in\scrX^{\circ}$ we have $\NC_{\bar{\scrC}}(x)=\{0\}$ and thus $g\in \scrA_{0}^{\bot}$. Hence, the claim follows from (c) of Proposition \ref{prop:propertiesV}.
\end{Proof}

%% file: PotentialReduction.tex
%%% Potential Reduction algorithm
%--------------------------------------------------------------------
% !TEX root = ./Main.tex
%

Based on the family of search directions $V(x,g)$, we now tailor the gradient input $g\in\Rn$ to derive a potential reduction algorithm solving problem \eqref{eq:P}. Throughout this paper we will work with a pair of functions $(f,h)$ such that: 
\begin{itemize}
\item[(i)] $h\in\scrF_{M,\nu}(\scrC)$ with some parameters $M>0$ and $\nu\in(2,4]$; 
\item[(ii)] $\cl\left(\dom h\right)=\bar{\scrC}$; 
\item[(iii)] $f:\Rn\to(-\infty,\infty]$ obeys Assumption \ref{ass:fcontinuous} and $\dom(h)\subseteq\dom(f)$. 
\end{itemize}
The next definition, due to \cite{BauBolTeb16,LuFreNes18} and \cite{BolSabTebVai18}, is fundamental to our analysis.

\begin{definition}
\label{def:L-adapt}
The pair of functions $(f,h)$ is \emph{L-smooth} if there exists a constant $L>0$ such that
\begin{equation}\label{eq:L-smooth}\tag{L}
f(y)-f(x)-\inner{\nabla f(x),y-x}\leq L D_{h}(y,x)\qquad\forall x,y\in \scrC 
\end{equation}
\end{definition}

It easy to check that $(f,h)$ being $L$-smooth is equivalent to $Lh-f$ being convex. Define the potential function 
\begin{equation}\label{eq:potential}
F_{\mu}(x)=f(x)+\mu h(x)\qquad\forall x\in\dom h.
\end{equation}

\begin{figure}
\centering
\includegraphics[width=0.5\textwidth]{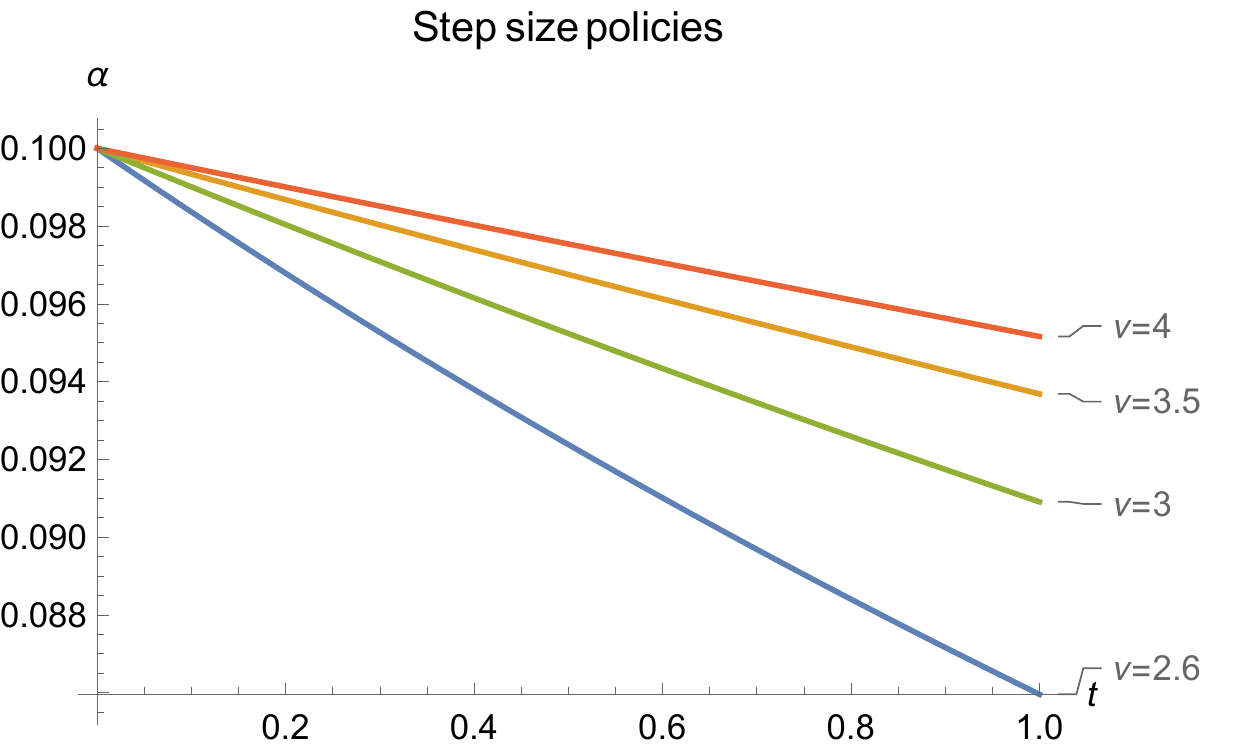}
\caption{Step size $\alpha=s_{\nu,L}(t)$ for $L+\mu=10$ and $\nu\in \{2.6,3,3.5,4\}$.}
\label{fig:step}
\end{figure}

If $(f,h)$ is an $L$-smooth pair, then the function $(L+\mu)h-F_{\mu}$ must be convex. Therefore, 
\begin{equation}\label{eq:descent}
F_{\mu}(x)\leq F_{\mu}(y)+\inner{\nabla F_{\mu}(y),x-y}+(L+\mu)D_{h}(x,y)\qquad\forall (x,y)\in\dom h\times\dom h.
\end{equation}
This inequality is in fact a modified descent lemma, in the spirit of \cite{BauBolTeb16}, for the non-convex, non-smooth composite function $F_{\mu}$. 
%%%%%%
\begin{remark}
If $\scrC=\Rn$ and $h(x)=\frac{1}{2}\norm{x}^{2}$, the pair $(f,h)$ is $L$-smooth if and only if the classical descent inequality 
\begin{align*}
f(x)-f(y)-\inner{\nabla f(y),x-y}\leq\frac{L}{2}\norm{x-y}^{2}
\end{align*}
holds for all $x,y\in\Rn$, i.e. the parameter $L$ is a surrogate for the Lipschitz constant of the Euclidean gradient map $x\mapsto\nabla f(x)$. 
\end{remark}
Define the search direction
\begin{equation}
d_{\mu}(x):=V(x,\nabla F_{\mu}(x))\qquad\forall x\in\scrX^{\circ}.
\end{equation}
From the first-order optimality condition of the search direction \eqref{eq:finder}, we know that 
\begin{equation}\label{eq:improve}
    \inner{\nabla F_{\mu}(x),d_{\mu}(x)}=-\norm{d_{\mu}(x)}^{2}_{x}\qquad\forall x\in\scrX^{\circ}.
\end{equation}
The associated dual variable is obtained by the evaluation of \eqref{eq:dual} as $y_{\mu}(x):=y(x,\nabla F_{\mu}(x))$. 
%%%%%%%%%%%%%%%%%%%%%%%%%%%%%%%%%%%%%
\begin{lemma}\label{lem:dual}
The dual function $y_{\mu}:\scrC\to\R^m$ is continuous.
\end{lemma}
\begin{Proof}
By \eqref{eq:dual}, the function $y_{\mu}(x)$ has the explicit expression $y_{\mu}(x)=(A[H(x)]^{-1}A^{\top})^{-1}A[H(x)]^{-1}\nabla F_{\mu}(x).$ Since $h\in\bC^{3}(\scrC)$, the matrix-valued mapping $H(x)=\nabla^{2}h(x)$ is continuous, and $x\mapsto\nabla F_{\mu}(x)$ is continuous as well on $\scrC$. The claim follows.
\end{Proof}
%%%%%%%%%%%%%%%%%%%%%%%%%
Define 
\begin{equation}
\lambda_{\mu}(x):=\norm{d_{\mu}(x)}_{x},\text{ and }\beta_{\mu}(x):=\norm{d_{\mu}(x)}_{2},
\end{equation}
and the transfer function 
%%%%%%%%%
\begin{equation}\label{eq:transfer}
T_{\mu}(x,\alpha):=x+\alpha d_{\mu}(x)\quad x\in\scrC,\alpha>0,\mu>0.
\end{equation}
%%%%%%%%%%%%
This mapping $T_{\mu}:\scrC\times[0,\infty)\to\Rn$ will serve as the generator of the numerical algorithm. Thanks to generalized self-concordance, we can easily determine the step length values $\alpha>0$ guaranteeing that $T_{\mu}(x,\alpha)\in \scrX^{\circ}$ for $x\in\scrX^{\circ}$. Indeed, by Lemma \ref{lem:Dikin}, a sufficient condition ensuring that we stay in the interior of the feasible set is to set $\alpha>0$ such that $\metric_{\nu}(x,T_{\mu}(x,\alpha))<1$. This leads to the bound $\alpha M\frac{\nu-2}{2}\lambda_{\mu}(x)^{\nu-2}\beta_{\mu}(x)^{3-\nu}<1.$ Let us define 
\begin{equation}\label{eq:delta}
\delta_{\mu}(x):= M\frac{\nu-2}{2}\lambda_{\mu}(x)^{\nu-2}\beta_{\mu}(x)^{3-\nu},
\end{equation}
so that $\metric_{\nu}(x,T_{\mu}(x,\alpha))=\alpha\delta_{\mu}(x)$ for all $x\in\scrX^{\circ},\alpha\geq 0$. Furthermore, if $\delta_{\mu}(x)>0$, we define $\bar{\alpha}_{\mu}(x):=1/\delta_{\mu}(x)$. Therefore, any choice of step size $\alpha\in(0,\bar{\alpha}_{\mu}(x))$, delivers a feasible step. Furthermore, for all $\alpha\in(0,\bar{\alpha}_{\mu}(x))$, we can apply the general descent inequality \eqref{eq:descent} to the $L$-smooth pair $(f,h)$, so that 
%\begin{align*}
%f(T_{\mu}(x,\alpha))\leq f(x)+\alpha \inner{\nabla f(x),d_{\mu}(x)}+LD_{h}(T_{\mu}(x,\alpha),x)\qquad\forall \alpha\in(0,\bar{\alpha}_{\mu}(x)).
%\end{align*} 
\eqref{eq:descent}, \eqref{eq:improve} and \eqref{eq:delta} give us the per-iteration estimate
\begin{align*}
F_{\mu}(T_{\mu}(x,\alpha))&\leq F_{\mu}(x)+\alpha \inner{\nabla F_{\mu}(x),d_{\mu}(x)}+(L+\mu)D_{h}(T_{\mu}(x,\alpha),x)\\
&=F_{\mu}(x)-\alpha \lambda_{\mu}(x)^{2}+(L+\mu)D_{h}(T_{\mu}(x,\alpha),x).
\end{align*}
Combining this with Lemma \ref{lem:Bregman}, we see that for all $x\in \scrC$ and $\alpha\in(0,\bar{\alpha}_{\mu}(x))$ 
\begin{align}
F_{\mu}(T_{\mu}(x,\alpha))&\leq F_{\mu}(x)-\alpha\lambda_{\mu}(x)^{2}+(L+\mu)\omega_{\nu}(\alpha\delta_{\mu}(x))\alpha^{2}\lambda_{\mu}(x)^{2}\nonumber\\
&=F_{\mu}(x)-\eta_{\mu}(x,\alpha),
\label{eq:progress}
\end{align}
where we have set 
\begin{equation}\label{eq:eta}
\eta_{\mu}(x,t):=t\lambda_{\mu}(x)^{2}-(L+\mu)\omega_{\nu}(t\delta_{\mu}(x))t^{2}\lambda_{\mu}(x)^{2}.
\end{equation}
Note that the barrier-generating kernel $h\in\scrF_{M,\nu}(\scrC)$ only appears in this per-iteration bound via the local norm of the search direction $\lambda_{\mu}(x)$. As such, the above bound can be seen as worst-case bound on the potential function decrease. This worst-case point of view is however very useful in determining an explicit step-size policy, akin to the recently proposed prox-based algorithms for convex composite self-concordant minimization \cite{CevKyrTra15,SunTran18}.
%%%%%%
\begin{proposition}\label{prop:step}
For all $x\in\scrX^{\circ},\mu,L>0$ and $\alpha\in(0,\bar{\alpha}_{\mu}(x))$, we have $T_{\mu}(x,\alpha)\in\scrX^{\circ}$. The optimal step-size rule, in the analytical worst-case sense, is given by 
\begin{equation}
\label{eq:alpha}
\alpha_{\mu}(x,L):=\left\{\begin{array}{ll} 
\frac{1}{\delta_{\mu}(x)}\left[1-\left(1+\frac{\delta_{\mu}(x)}{L+\mu}\frac{4-\nu}{\nu-2}\right)^{-\frac{\nu-2}{4-\nu}}\right] & \text{ if }\nu\in(2,3)\cup(3,4),\\
\frac{1}{\delta_{\mu}(x)+L+\mu} & \text{ if }\nu=3,\\
\frac{1}{\delta_{\mu}(x)}\left[1-\exp\left(-\frac{\delta_{\mu}(x)}{L+\mu}\right)\right] & \text{if }\nu=4. 
\end{array}\right. 
\end{equation}
\end{proposition}
%%%%%%%%%%%%%%%%%%%%%%%%%%%%%
The proof of this Proposition is a rather technical computation, and therefore delegated to Appendix \ref{app:Step}. It is however interesting to note that the self-concordance parameter $\nu$ plays a somewhat symmetric role around its values $\nu\in(2,3)\cup(3,4)$. Moreover, it is remarkable that the theoretical upper bound on the step size ensuring feasibility, $\bar{\alpha}_{\mu}(x)$, is independent of the constant $L$. It appears only when we compute the optimal step size $\alpha_{\mu}(x,L)$, which, in turn, is available in a closed-form expression. This functional form of the step size policy allows for a direct comparison in dependence of the generalized self-concordance parameter $\nu\in(2,4]$. In Figure \ref{fig:step} we provide a numerical illustration on the ordering of the step sizes, realizing for fixed parameter pair $(\mu,L)$, we can think of the function $\alpha_{\mu}(x,L)$ as the output of a function of the composition $s_{\mu,L,\nu}\circ \delta_{\mu}$, where $s_{\mu,L,\nu}:(0,\infty)\to(0,\infty)$ is given by 
\[
s_{\mu,L,\nu}(t):=\left\{\begin{array}{ll} 
\frac{1}{t}\left[1-\left(1+\frac{t}{L+\mu}\frac{4-\nu}{\nu-2}\right)^{-\frac{\nu-2}{4-\nu}}\right] & \text{ if }\nu\in(2,3)\cup(3,4),\\
\frac{1}{t+L+\mu} & \text{ if }\nu=3,\\
\frac{1}{t}\left[1-\exp\left(-\frac{t}{L+\mu}\right)\right] & \text{if }\nu=4. 
\end{array}\right. 
\]

We see that larger parameters $\nu$ lead to higher step sizes and thus to more aggressive schemes. This clearly indicates that the choice of the metric-inducing kernel matters in the design of the algorithm. Observe that the function $s_{\mu,L,\nu}$ is well defined at $t=0$ and attains the same value $1/(L+\mu)$ for all $\nu\in(2,4]$. 

Both parameters $\mu$ and $L$ are seen to have the same effect on the step size policy: Larger values imply smaller step sizes. Hence, for optimization purposes it is of utmost interest to pick these parameters in a way that prevents $\HBA(\mu,L)$ making to small steps. However, the two parameters also play different roles in the design of the method. While $\mu$ is a \emph{barrier parameter} guaranteeing that the algorithm is an interior-point method, the parameter $L$ is dictated by the pair $(f,h)$ in order to guarantee the descent property \eqref{eq:L-smooth}. Intuitively, we would like to run the algorithm with a numerically small value of $\mu$. The descent property \eqref{eq:L-smooth} tells us that for any $\tilde{L}\geq L$, we can guarantee a sufficient decrease in the potential function, so the design question becomes how small the the parameter $L$ can be chosen with a guarantee to obtain a sufficient decrease. In order to answer these questions, we will develop an adaptive version of the base algorithm $\HBA(\mu,L)$. Since the analysis of this adaptive version will rely on general results obtained for the base scheme $\HBA(\mu,L)$, we start our mathematical analysis with the assumption that both parameters are fixed, and later make their choice dynamic. 

%% file: algo.tex
%----------------------------------------------------------------------
%%% Algorithm
%----------------------------------------------------------------------
% !TEX root = ./Main.tex
%

Let $\const>0$ be a positive constant. For the construction of our algorithmic scheme it will be important to have access to a $\const$-analytic center, i.e. a point $x^{0}\in\scrX^{\circ}$ such that 
\begin{equation}\label{eq:h}
h(x)\geq h(x^{0})-\const\qquad\forall x\in\scrX. 
\end{equation}
To obtain such a point $x^{0}$, we need to be able to approximately solve the minimization problem 
\[
\min_{x\in\scrX}h(x). 
\]
This assumption is very common in potential reduction schemes \cite{Ye92,Ye98,HaeLiuYe18}. In case where $\scrX$ is bounded, existence of an exact analytic center is guaranteed and we can use any efficient solver for computing it. In the case where $\scrX$ is unbounded, existence of an exact analytic center is guaranteed if the classical existence condition given by the Weierstrass' theorem are satisfied. If $h\in\scrF_{M,\nu}(C)$ for $\nu\in (3,4]$ then by Definition \ref{def:bgk} the function $h$ is coercive, and therefore the program \eqref{eq:h} always has a solution. In the remaining cases, where $\nu\in(2,3]$ we either have to assume that the feasible set is bounded or that there exists a point $x\in\scrX$ under which the dual norm $\norm{\nabla h(x)}^{\ast}_{x}$ is "small". A precise meaning of this statement can be given by adapting the relevant arguments in \cite{SunTran18} to the current geometry with linear equality constraints and leave this, somewhat off-topic, exercise to the reader.

\subsection{HBA iterations for fixed $\mu$}
\label{sec:fixed}
Given the general template described in Section \ref{sec:general}, the first algorithmic scheme of interest in this paper is easy to describe. Assuming the standing hypothesis Assumptions \ref{ass:C}-\ref{ass:fcontinuous} in place, we are given a pair of functions $(f,h)$ such that $h\in\scrF_{M,\nu}(\scrC)$, and $L$-smoothness holds for some $L\geq 0$. Given a $\const$-approximate analytic center $x^{0}\in\scrX^{\circ}$ as initial condition, we generate a sequence $(x^{k})_{k\geq 0}$ recursively by 
\[
x^{k+1}=T_{\mu}(x^{k},\alpha_{k})=x^{k}+\alpha_{k}v^{k}_{\mu}\qquad\forall k\geq 0,
\]
where $T_{\mu}(x,\alpha)$ is the transfer function defined in \eqref{eq:transfer}, with step-size $\alpha_{k}=\alpha_{\mu}(x^{k},L)$ and search direction $v^{k}_{\mu}=d_{\mu}(x^{k})$.
%\begin{equation}\label{eq:step_fixed}
%\alpha_{k}=\alpha_{\mu}(x^{k},L)\qquad\forall k\geq 0,
%\end{equation}
%and search direction
%\begin{equation}\label{eq:search_fixed}
%v^{k}_{\mu}=d_{\mu}(x^{k})\qquad \forall k\geq 0,
%\end{equation}
%and a step-size 
%\begin{equation}\label{eq:step_fixed}
%\alpha_{k}=\alpha_{\mu}(x^{k},L)\qquad\forall k\geq 0.
%\end{equation}
%The resulting iterates of $\HBA(\mu,L)$ read therefore as the discrete-time dynamical system
%\begin{align*}
%&x^{k+1}=T_{\mu}(x^{k},\alpha_{k})\qquad\forall k\geq 0,\\
%&x^{0}\in\scrX^{0}\text{ given.}
%\end{align*}
The pseudo-code corresponding to the conceptual implementation of $\HBA(\mu,L)$ reads as Algorithm \ref{alg:HBA}.
%%%%%%%%%%%%%%%%%%%%%%%%%%%%%%%%%%%%%%%%%%%%%%%%%%%%%%%%%%%
%%%%%%%%%%%%% Pseudo Code of HBA%%%%%%%%%%%%%%%%%%%%%%%%%%%%%%%%%%%%%%
%%%%%%%%%%%%%%%%%%%%%%%%%%%%%%%%%%%%%%%%%%%%%%%%%%%%%%%%%%%%
\input{HBA_fixedmu}
%%%%%%%%%%%%%%%%%%%%%%%%%%
Let us make some remarks on the computational efficiency of $\HBA(\mu,L)$. Since $\scrC$ contains no lines and the sequence $(x^{k})_{k\geq 0}$ stays in the relative interior for the optimization problem's feasible set, we guarantee that $H(x^{k})\succ 0$ for every iteration of the algorithm. Hence, the main computational step in $\HBA(\mu,L)$ \eqref{eq:finder} is always well-posed while executing the algorithm, and delivers a unique solution. The complexity of \eqref{eq:finder} is the same as for a Newton method; It requires $\bigoh(n^3)$ operations via either a Cholesky decomposition, or a well implemented conjugate gradient (CG) method. Hence, in terms of per-iteration complexity $\HBA(\mu,L)$ is comparable with Newton methods. In many applications, however, the matrix $H(x^{k})$ has a special structure which makes the application of heavy analytic machinery unnecessary. Indeed, most barrier functions $h$ used in the literature are additively separable and the resulting Hessian matrix $H$ is therefore diagonal. In this case, the computational complexity is essentially determined by the density of the matrix $A$, and in many applications (e.g. resource allocations problems where $A$ embodies network flow constraints) we will be able to implement a closed-form expression for the search direction $v^{k}_{\mu}$. Hence, in such favorable instances, the per-iteration computational overhead of implementing $\HBA(\mu,L)$ is rather small.

\subsection{Adaptive HBA}
\label{sec:adaptive}
The basic algorithmic scheme $\HBA(\mu,L)$ utilizes knowledge of the barrier parameter $\mu$ and the L-smoothness parameter $L>0$ in order to determine the step size $\alpha_{\mu}(x,L)$. Knowing a-priori the parameter $\mu$ is not a very demanding, since it is chosen by the user at the beginning of the implementation. However, running the basic \acl{HBA} scheme with a stiff parameter $L$ might be inefficient since it forces us to rescale the step size with the same constant factor $\mu+L$ globally. This might lead to unnecessary small steps, resulting in long run times of the method. To overcome this drawback we present in this section a new adaptive method of \acl{HBA}, where all necessary information about $L$ can be accumulated by an appropriate "line-search" strategy. The thus resulting \acl{AHBA} method closely resembles ideas spelled out in "universal gradient methods" defined by Nesterov in \cite{Nes15}. 
%%%%%%%%%%%%%%%%%%%%%%%%%%%%
%%%%%%%%%%% Adaptive Version HBA%%%%%%%
%%%%%%%%%%%%%%%%%%%%%
\input{AHBA}
%%%%%%%%%%%%%%%%%%%%%%%%%%%%%%%%%%%%%%%%%%%%%%%
%%%%%%%%%%%%%%%%%%%%%%%%%%%%%%%%%%%%

%The intuition behind the design of $\AHBA(\mu,L)$ is the following. Recall that a basic compatibility assumption of our method is the $L$-smoothness of the pair $(f,h)$. This condition ensures the descent inequality \eqref{eq:L-smooth} uniformly over the set $\scrX^{\circ}$. The basic \ac{HBA} mechanism fixes the search direction $d_{\mu}(x)$, and the only remaining question is to sample a point along the ray $\{T_{\mu}(x,\alpha),0\leq\alpha<\bar{\alpha}_{\mu}(x)\}$ to guarantee significant progress of the method. Implementing the optimal step-size policy $\alpha_{\mu}(x,L)$ in this restricted line search leaves us with the choice of the parameter $L$. Since the step size is inversely related to the Lipschitz-like parameter $L$, it would be beneficial to overall efficiency the algorithm if we are able to use a smaller parameter $L'$ than the uniform value $L$. This is exactly what we search for in the subroutine \eqref{eq:line}.

%% file: HBA_fixedmu.tex
%----------------------------------------------------------------------
%%% HBA-fixed mu
%----------------------------------------------------------------------
% !TEX root = ../Main.tex
%

\begin{algorithm}[t]
\caption{$\HBA(\mu,L)$}
\label{alg:HBA}
\SetAlgoLined
\KwData{
kernel generating distance $h\in\scrF_{M,\nu}(\scrC)$ such that $(f,h)$ is $L$-smooth\; 
Barrier parameter $\mu>0$.
}
\KwResult{Stationary point of $F_{\mu}$.}
Initial point: $\const$-analytic center $x^{0}\in\scrX^{\circ}$\;
 \While{$k=0,1,\ldots,k_{\max}$}{
  obtain $x^{k}$ and $\nabla F_{\mu}(x^{k})$\; 
    \eIf{Stopping condition not satisfied}{
Solve the linear system \eqref{eq:finder} for $x=x^{k}$ and $g=\nabla F_{\mu}(x^{k})$\;
Denote by $(v_{\mu}^{k},y_{\mu}^{k})^{\top}$ the solution\;
Compute step size $\alpha_{k}=\alpha_{\mu}(x^{k},L)$\; 
Update $x^{k+1}=x^{k}+\alpha_{k}v_{\mu}^{k}$.
}{
   Stop and report $x^{k}$ as the solution
  }
 }
\end{algorithm}

%% file: AHBA.tex
%----------------------------------------------------------------------
%%% HBA-adaptive
%----------------------------------------------------------------------
% !TEX root = ../Main.tex
%

\begin{algorithm}[t]
\caption{$\AHBA(\mu)$}
\label{alg:AHBA}
\SetAlgoLined
\KwData{
kernel generating distance $h\in\scrF_{M,\nu}(\scrC)$ such that $(f,h)$ is $L$-smooth\; 
Barrier parameter $\mu>0$.
}
\KwResult{Stationary point of $F_{\mu}$.}
Initial point: $\const$-analytic center $x^{0}\in\scrX^{\circ}$\;
 \While{$k=0,1,\ldots,k_{\max}$}{
  obtain $x^{k}$\; 
    \eIf{Stopping condition not satisfied}{
   Solve the linear system  \eqref{eq:finder} for $x=x^{k}$ and $g=\nabla F_{\mu}(x^{k})$.\;
Denote by $(v_{\mu}^{k},y_{\mu}^{k})^{\top}$ the solution\;
 Find the smallest $i_k \geq 0$ such that $z^{k}=x^{k}+\alpha_{\mu}(x^{k},2^{i_k-1}L_{k})v_{\mu}^{k}$ satisfies
	\begin{equation}
		\label{eq:line}
		f(z^{k}) \leq f(x^k) + \inner{\nabla f(x^{k}),x^{k+1}-x^{k}}+2^{i_k-1}L_{k}D_{h}(z^{k},x^{k}).
	\end{equation}
Set $L_{k+1} = 2^{i_k-1}L_{k}$\;
Update $x^{k+1}=z^{k}$.
}{
   Stop and report $x^{k}$ as the solution.
  }
 }
\end{algorithm}

%% file: convergence.tex
%%% Convergence analysis
%--------------------------------------------------------------------
% !TEX root = ./Main.tex
%

We organize our discussion on the long-run properties of $\HBA(\mu,L)$ and $\AHBA(\mu)$ in two parts. The first part is concerned with the asymptotic convergence properties of the two methods. In the second part we will discuss the non-asymptotic complexity properties of the method. Throughout this section we assume that Assumptions \ref{ass:C}-\ref{ass:fcontinuous} are satisfied. 

\subsection{Asymptotic convergence}
\label{sec:asym}
Let $(x^{k})_{k\geq 0}$ be a sequence generated by $\HBA(\mu,L)$, with search direction $v^{k}_{\mu}=d_{\mu}(x^{k})$ and step-size policy $\alpha_{k}=\alpha_{\mu}(x^{k},L)$. Let us introduce the associated sequences $(\lambda_{k})_{k\geq 0},(\beta_{k})_{k\geq 0},(\delta_{k})_{k\geq 0}\subset[0,\infty)$ by 
\begin{align}
\lambda_{k}&:=\norm{v_{\mu}(x^{k})}_{x^{k}},\; \beta_{k}:=\norm{v_{\mu}(x^{k})}_{2},\; \delta_{k}:=M\left(\frac{\nu-2}{2}\right)\lambda^{\nu-2}_{k}\beta^{3-\nu}_{k}. \label{eq:lambda_beta_delta}
\end{align}
Let us define the per-iteration progress along the thus produced sequence, quantified in \eqref{eq:progress}, as 
\begin{equation}\label{eq:Delta}
\Delta_{k}:=\eta_{\mu}(x^{k},\alpha_{k})\qquad\forall k\geq 0.
\end{equation}
From Proposition \ref{prop:step}, we immediately deduce an online version of the step-sizes together with a descent inequality for the potential function, summarized in the next Proposition.
%%%%%%%%%%%%%%%%%%%%%%%%%%%%%%%%%%%%%
\begin{proposition}\label{prop:descent_fixed}
Let $(x^{k})_{k\geq 0}$ be the sequence generated by $\HBA(\mu,L)$ with the step size policy 
\begin{equation}
\label{eq:alpha_k}
\alpha_{k}:=\left\{\begin{array}{ll} 
\frac{1}{\delta_{k}}[1-\left(1+\frac{\delta_{k}}{L+\mu}\frac{4-\nu}{\nu-2}\right)^{-\frac{\nu-2}{4-\nu}}] & \text{ if }\nu\in(2,3)\cup(3,4),\\
\frac{1}{\delta_{k}+L+\mu} & \text{ if }\nu=3,\\
\frac{1}{\delta_{k}}\left[1-\exp\left(-\frac{\delta_{k}}{L+\mu}\right)\right] & \text{if }\nu=4. 
\end{array}\right. 
\end{equation}
Then, for all $k\geq 0$ we have
\begin{equation}\label{eq:decrease}
F_{\mu}(x^{k+1})\leq F_{\mu}(x^{k})-\Delta_{k},
\end{equation}
where $\Delta_{k}$ is defined in \eqref{eq:Delta}. Moreover, this step size rule is optimal in the worst-case analytic sense.
\end{proposition}
%%%%%%%%%%%%%%%%%%%%%%%%
We next provide some general properties of $\HBA(\mu,L)$.
%%%%%%%%%%%%%%%%%%%%%%%%%%%%%%%%%%
\begin{proposition}
\label{prop:Delta}
Let $(x^{k})_{k\geq 0}$ be generated by $\HBA(\mu,L)$, and set $f^{\ast}:=\inf_{x\in\scrX}f(x)$. Then, the following assertions hold:
\begin{itemize}
\item[(a)] $\left(F_{\mu}(x^{k})\right)_{k\geq 0}$ is non-increasing;
\item[(b)] $\sum_{k\geq 0}\Delta_{k}<\infty$, and hence the sequence $\left(\Delta_{k}\right)_{k\geq 0}$ converges to 0;
\item[(c)] $\min_{0\leq k<K}\Delta_{k}\leq\frac{1}{K}[f(x^{0})-f^{\ast}+\mu\const]$.
\end{itemize}
\end{proposition}
\begin{Proof}
Unraveling the expressions in eq. \eqref{eq:descent}, we get for all $k\geq 0$,
\begin{align*}
f(x^{k+1})-f(x^{k})\leq -\Delta_{k}+\mu[h(x^{k})-h(x^{k+1})]. 
\end{align*}
Telescoping this expression shows that for all $K\geq 1$, 
\begin{align*}
f(x^{K})-f(x^{0})\leq -\sum_{k=0}^{K-1}\Delta_{k}+\mu[h(x^{0})-h(x^{K})]. 
\end{align*}
Since $x^{0}$ is a $\const$-analytic center, the left-hand side in the above display can be majorized to obtain the bound
\begin{align*}
f(x^{K})-f(x^{0})\leq-\sum_{k=0}^{K-1}\Delta_{k}+\mu\const. 
\end{align*}
Since $\Delta_{k}>0$ , the sequence $\left(F_{\mu}(x^{k})\right)_{k\geq0}$ is monotonically decreasing. Since $h(x)-h(x^{0})\geq -\const$, and $f$ is bounded from below, the potential function $F_{\mu}$ is bounded from below as well. Therefore $\liminf_{k\to\infty}F_{\mu}(x^{k})=\lim_{k\to\infty}F_{\mu}(x^{k})$ exists and equals a number $F^{\ast}_{\mu}\in(-\infty,\infty)$. It follows that $(x^{k})_{k\geq 0}\subset\dom h=\scrC$, and therefore $\lim_{k\to\infty} f(x^{k})$ exists as well. Calling $f^{\ast}:=\inf\{f(x): x\in\scrX\}>-\infty$, we conclude that for all $K\geq 1$, 
\begin{equation}
\sum_{k=0}^{K-1}\Delta_{k}\leq f(x^{0})-f(x^{K})+\mu \const\leq f(x^{0})-f^{\ast}+\mu\const,
\end{equation}
and 
\begin{equation}
\min_{1\leq k\leq K}\Delta_{k}\leq\frac{1}{K}[f(x^{0})-f^{\ast}+\mu\const].
\end{equation}
Hence, $\lim_{k\to\infty}\Delta_{k}=0$. 
\end{Proof}
%%%%%%%%%%%%%%%%%%%%%%
We turn now to the convergence properties of $\HBA(\mu,L)$. Our aim is to show that accumulation points of the sequence $(x^{k})_{k\geq 0}$ generated by the algorithm are stationary points of the potential function $F_{\mu}$. We start by proving some auxiliary results.
%%%%%%%%%%%%%%
\begin{lemma}\label{lem:bounded}
Let $(x^{k})_{k\geq 0}$ be generated by $\HBA(\mu,L)$. Then, $(x^{k})_{k\geq 0}$ is bounded. 
\end{lemma}
\begin{Proof}
Since $\left(F_{\mu}(x^{k})\right)_{k\geq 0}$ is monotonically decreasing, we have 
\begin{align*}
f(x^{k+1})-f(x^{k})&\leq -\Delta_{k}-\mu[h(x^{k+1})-h(x^{k})].
\end{align*}
Hence, for all $K\geq 1$, using the $\const$-analytic center property of the initial condition $x^{0}$, we get
\[
f(x^{K})\leq f(x^{0})+\mu\const. 
\]
Hence, $x^{k}\in\lev_{f}(f(x^{0})+\mu\const)$. Since $f$ has bounded level sets (Assumption \ref{ass:levelbound}), the entire sequence $(x^{k})_{k\geq 0}$ is bounded. 
\end{Proof}
%%%%%%%%%%%%
Define the limit set 
\begin{equation}
\omega(x^{0}):=\{p\in\scrX: \exists (k_{q})_{q\in\N}\uparrow\infty,\lim_{k_{q}\to\infty}x^{k_{q}}=p\}. 
\end{equation}
Thanks to Lemma \ref{lem:bounded}, standard results imply that $\omega(x^{0})$ is nonempty, connected and compact (see e.g. \cite[Lem.5]{BolShaTab14}). Furthermore, $\lim_{k\to\infty}\dist(x^{k},\omega(x^{0}))=0$. For $x\in\scrX^{\circ}$, define 
\begin{equation}
\scrS_{\mu}(x):=\{y\in\scrX: F_{\mu}(y)\leq F_{\mu}(x)\}=\lev_{F_{\mu}}(F_{\mu}(x))\cap\scrA.
\end{equation}
Since $\HBA(\mu,L)$ is a descent method for the potential function $F_{\mu}$, we immediately conclude that $(x^{k})_{k\geq 0}\subseteq\scrS_{\mu}(x^{0})$. 
%%%%%%%%%%%%%%%%%%%%%%%
\begin{lemma}\label{lem:S_compact}
Let $x^{0}$ be a $\const$-analytic center. Then $\scrS_{\mu}(x^{0})$ is a compact subset in $\scrX^{\circ}=\scrC\cap\scrA$.
\end{lemma}
\begin{Proof}
Note that
\begin{align*}
\scrS_{\mu}(x^{0})&=\{x\in\scrX: f(x)\leq f(x^{0})+\mu[h(x^{0})-h(x)]\}\\
&\subseteq\{x\in\scrX: f(x)\leq f(x^{0})+\mu\const\}\\
&=\lev_{f}(f(x^0)+\mu\const)\cap\scrA. 
\end{align*}
Since $f$ has bounded level sets (Assumption \ref{ass:levelbound}), the set $\scrS_{\mu}(x^{0})$ is bounded as well. It remains to prove that the set $\scrS_{\mu}(x^{0})$ is closed. To that end, let $(x_{j})_{j\geq 1}$ be a converging sequence with $\lim_{j\to\infty}x_{j}=\bar{x}$ and $x_{j}\in\scrS_{\mu}(x^{0})$ for all $j\geq 1$. Then, $f(x^{j})+\mu h(x^{j})\leq r\equiv f(x^{0})+\mu h(x^{0})$ for all $j\geq 1$. If $\bar{x}\in\bd(\scrC)$, then $h(x_{j})\to\infty$ and we immediately obtain a contradiction. Hence, $x\in\scrC\cap\scrA$, and the restriction of the composite function $f+\mu h$ on this domain is continuous. We conclude $\bar{x}\in\scrS_{\mu}(x^{0})$.
\end{Proof}

\begin{corollary}\label{cor:interior}
$\omega(x^{0})\subset\scrX^{\circ}$.
\end{corollary}

Let $\sigma_{\min}(x)$ denote the smallest and $\sigma_{\max}(x)$ the largest eigenvalue of the Hessian $H(x)=\nabla^{2}h(x)$. Since $H(x)\succ 0 $ for all $x\in\scrX^{\circ}$, we conclude that $\sigma_{\min}(x)>0$. Moreover, the compactness of the set $\scrS_{\mu}(x^{0})$ allows us to define the positive constant
\begin{equation}\label{eq:sigma}
\sigma_{h}:=\min_{x\in\scrS_{\mu}(x^{0})}\sigma_{\min}(x)
\end{equation}
Hence, along the iterates of $\HBA(\mu,L)$, we have 
\begin{equation}\label{eq:norm_equal}
\lambda_{k}\geq \sqrt{\sigma_{h}}\beta_{k}\qquad\forall k\geq 0.
\end{equation}
Since
\[
\alpha_{k}\delta_{k}=\metric_{\nu}(x^{k},x^{k+1})=M(\frac{\nu}{2}-1)\lambda_{k}^{\nu-2}\beta_{k}^{3-\nu},
\]
the following lower and upper bounds can be established for $\nu\in(2,3]$:
\begin{equation}\label{eq:sandwich}
M\left(\frac{\nu}{2}-1\right)\sigma_{h}^{\frac{\nu-2}{2}}\alpha_{k}\beta_{k}\leq\alpha_{k}\delta_{k}\leq M\left(\frac{\nu}{2}-1\right)\sigma_{h}^{-\frac{3-\nu}{2}}\alpha_{k}\lambda_{k}.
\end{equation}
This inequality will be key to prove convergence of the method to a stationary point of the potential function when $\nu\in(2,3]$. For $\nu\in(3,4]$, we will need to upper bound the local norm of the search direction, $\lambda_{k}$, as well. Let $\sigma_{\max}(x)\in(0,\infty]$ be the largest eigenvalue of the Hessian matrix $H(x)$. Since $\scrS_{\mu}(x^{0})$ is a compact set in $\scrX^{\circ}$, the quantity 
\begin{equation}\label{eq:tau}
\tau_{h}:=\max_{x\in\scrS_{\mu}(x^{0})}\sigma_{\max}(x)
\end{equation}
is well-defined and finite. Given these bounds, we see that for all $x\in\scrS_{\mu}(x^{0})$ we have 
\begin{equation}
\sigma_{h}\Id\preceq H(x)\preceq\tau_{h}\Id
\end{equation}
so that the function $h$ is $\sigma_{h}$-smooth and $\tau_{h}$-strongly convex on the compact set $\scrS_{\mu}(x^{0})$. The quantity $\kappa_{h}=\frac{\tau_{h}}{\sigma_{h}}\geq 1$ is the condition number of $h$. Hence, along the sequence $(x^{k})_{k\geq 0}$ generated by $\HBA(\mu,L)$, we can upper bound the local norm of the search direction by
\begin{equation}\label{eq:norm_equal2}
\lambda_{k}\leq\sqrt{\tau_{h}}\beta_{k}\qquad\forall k\geq 0.
\end{equation}
All these estimates together will be needed to prove the main result of this section, represented by the following Theorem. 
%%%%%%%%%
\begin{theorem}\label{th:gradient}
Let $(x^{k})_{k\geq 0}$ be the sequence generated by $\HBA(\mu,L)$ with step-size policy $(\alpha_{k})_{k\geq 0}$ described in \eqref{eq:alpha_k}. Then, $\omega(x^{0})\subseteq \{x\in\scrX:(\exists y\in\R^{m}): \nabla F_{\mu}(x)-A^{\top}y=0\}$.   
\end{theorem}
\begin{Proof}
See Appendix \ref{app:stationary}.
\end{Proof}
%%%%%%
As a consequence of this Theorem, it follows that the trajectory $(x^{k})_{k\geq 0}$ exhibits a decaying energy in the metric-like function $\metric_{\nu}$: 
\begin{corollary}
$\lim_{k\to\infty}\metric_{\nu}(x^{k},x^{k+1})=0$. 
\end{corollary}
\begin{Proof}
By definition, $\metric_{\nu}(x^{k},x^{k+1})=\alpha_{k}\delta_{k}$ for all $k\geq 0$ and $\nu\in(2,4]$. In Appendix \ref{app:stationary} we have shown that $\liminf_{k\to\infty}\alpha_{k}>0$ and $\limsup_{k\to\infty}\delta_{k}=0$. The claim follows. 
\end{Proof}
%%%%%%%%%%%%%%%%%%%%
\subsection{Non-asymptotic bounds}
In this section we provide complexity estimates for the non-adaptive base algorithm $\HBA(\mu,L)$. To do so, we report first a useful technical corollary of the proof of Theorem \ref{th:gradient}. 
\begin{lemma}\label{lem:omega}
 Let $(x^{k})_{k\geq 0}$ be generated by $\HBA(\mu,L)$, with corresponding potential reduction sequence $(\Delta_{k})_{k\geq 0}$ defined in \eqref{eq:Delta}. For each generalized self-concordance parameter $\nu\in(2,4]$, there exists a strictly increasing function $\tilde{\omega}_{\nu}:(0,\infty)\to(0,\infty)$ satisfying 
\begin{equation}
\Delta_{k}\geq\tilde{\omega}_{\nu}(\lambda_{k})\qquad\forall k\geq 0.
\end{equation}
In particular, this function is given by 
\begin{align*}
\tilde{\omega}_{\nu}(t):&=\left\{\begin{array}{ll}
\tilde{\gamma}_{\nu}t\min\left\{\frac{2\sigma_{h}^{\frac{3-\nu}{2}}}{M(\nu-2)},\frac{t}{-\cb(L+\mu)}\right\} & \text{if }\nu\in(2,3),\\
\tilde{\gamma}_{\nu}t\min\left\{\frac{2}{M(\nu-2)}\tau_{h}^{-\frac{3-\nu}{2}},\frac{t}{-\cb(L+\mu)}\right\} & \text{if }\nu\in(3,4),\\
\frac{2(1-\ln(2))t}{M(L+\mu)}\min\left\{(L+\mu),\frac{M}{2}t\right\} & \text{if }\nu=3,\\
t\exp(-1)\min\left\{\frac{1}{\sqrt{\tau_{h}}M},\frac{t}{L+\mu}\right\} & \text{if }\nu=4,
\end{array}\right.
\end{align*}
where 
\begin{align*}
&\cb:=\frac{2-\nu}{4-\nu}\text{ for } \nu\in(2,4), \text{ and }\\
&\tilde{\gamma}_{\nu}:=1+\frac{4-\nu}{2(3-\nu)}\left(1-2^{\frac{2(3-\nu)}{4-\nu}}\right)\text{ for  }\nu\in(2,3)\cup(3,4).
\end{align*} 
\end{lemma}
\begin{Proof}
The proof follows from eqs. \eqref{eq:finalDelta1}, \eqref{eq:finalDelta2}, \eqref{eq:finalDelta3}, and \eqref{eq:finalDelta4} in Appendix \ref{app:stationary}.
\end{Proof}
A remarkable observation we can make from this Corollary is that the eigenvalue bounds number $\tau_{h}$ and $\sigma_{h}$ only appear for the generalized self-concordance parameters $\nu\in(2,4]\setminus\{3\}$.

Lemma \ref{lem:omega} is key to prove the first iteration complexity bounds to estimate the number of steps needed to ensure that the local norm of the search direction is smaller than a user-defined tolerance. In the context of proximal algorithms for solving composite self-concordant minimization problems with convex data, a similar result has been established by \cite{CevKyrTra15}. We instead derive such a basic complexity estimate in the setting of Hessian-barrier methods for non-convex optimization problems without Lipschitz gradient assumptions and generalized self-concordant penalties. 
%%%%%%%%%%%%%%%%%%%%%%%%%%%%%%%%%
\begin{lemma}\label{lem:complex_basic}
Suppose Assumptions \ref{ass:C}-\ref{ass:fcontinuous} hold. Let $(x^{k})_{k\geq 0}$ be the sequence generated by $\HBA(\mu,L)$. Define the stopping time 
\begin{equation}
\bN(\eps,x^{0},\nu,L):=\min\left\{k\geq 0:\lambda_{k}<\eps\right\}.
\end{equation}
Then, 
\begin{equation}
\bN(\eps,x^{0},\nu,L)\leq \left\lceil \frac{f(x^{0})-f^{\ast}+\mu\const}{\tilde{\omega}_{\nu}(\eps)}\right\rceil.
\end{equation}
\end{lemma}
\begin{Proof}
By definition, for all $0\leq k\leq\bN(\eps,x^{0},\nu,L)-1$,we have $\lambda_k \geq \eps$ and, due to the strong monotonicity of the function $\tilde{\omega}_{\nu}$, that $\tilde{\omega}_{\nu}(\lambda_{k})\geq\tilde{\omega}_{\nu}(\eps)$. Therefore, using the per-iteration descent of the potential function given by 
\[
F_{\mu}(x^{k+1})-F_{\mu}(x^{k})\leq-\Delta_{k}\leq-\tilde{\omega}_{\nu}(\lambda_{k})\leq - \tilde{\omega}_{\nu}(\eps),
\]
we readily conclude for $N > \bN(\eps,x^{0},\nu,L)$, 
\[
f^{\ast}\leq f(x^{N})\leq f(x^{0})-N\tilde{\omega}_{\nu}(\eps)+\mu\const< f^*. 
\]
%\PD{Thus, $N \leq \bN(\eps,x^{0},\nu,L)$.}
Solving for $N$ gives the claimed bound.
\end{Proof}
%%%%%%%%%%%
Our second iteration complexity result gives a more precise estimate on the number of steps needed to make the local norm of the search direction as small as desired. In particular, the next estimate provides us with an easy-to-implement stopping criterion for $\HBA(\mu,L)$, building on the insights gained from Lemma \ref{lem:complex_basic}. Let $\eps>0$ be a target precision level, specified before the algorithm is started, and set $\mu=4\eps$. We elect to terminate $\HBA(4\eps,L)$ whenever $F_{4\eps}(x^{K+1})-F_{4\eps}(x^{K})\geq -\frac{\hat{\gamma}_{\nu}\eps^{2}}{L+4\eps}$ %\PD{Maybe here it is better to define $\tilde{\omega}_{\nu}(t,\mu)$ and use stopping criterion $F_{4\eps}(x^{K+1})-F_{4\eps}(x^{K})\geq \tilde{\omega}_{\nu}(\eps,4\eps)$? Then by the monotonicity of $\tilde{\omega}_{\nu}(t,4\eps)$ at the stopping time, $\lambda_{K}\leq\eps$.} 
at iteration $K$, and report the iterate $x^{K}$. When this happens for the first time, we will show that $\lambda_{K}\leq\eps$. If this stopping criterion is not satisfied, we continue with the execution of the protocol $\HBA(4\eps,L)$ until an upper bound on the number of iterations $K=O(\eps^{-2})$ is reached.  Implementing this stopping criterion, we therefore are guaranteed to reach a point $x^{K}$ either satisfying $f(x^{K})-\inf_{x\in\scrX}f(x)= f(x^{K})-f^{\ast}\leq \eps$, or else $\lambda_{K}\leq \eps$. Together with this stopping criterion we see that $\HBA(4\eps,L)$ solves a constrained problem with potential non-differentiability at the boundary, with an iteration complexity of $O(\eps^{-2})$. For this type of problem, such a rate is the best known in the literature \cite{NesConvex}. It is also worth emphasizing that the transfer map $T_{\mu}$ admits a closed form expression and the complexity of computing this map is on the same order as standard Newton methods.

We are now ready to present the main complexity result for $\HBA(4\eps,L)$. Define the constant
\begin{equation}
\hat{\gamma}_{\nu}:=\left\{\begin{array}{cl} 
\tilde{\gamma}_{v}\frac{4-\nu}{\nu-2} & \text{if }\nu\in(2,3)\cup(3,4)\\
1 & \text{if }\nu=3,\\
\exp(-1) & \text{if }\nu=4.
\end{array}\right.
\end{equation}
Consider the stopping time 
\begin{equation}\label{eq:K1}
\bK_{1}(\eps,x^{0},\nu,L):=\min\left\{k\geq 0: \tilde{\omega}_{\nu}(\lambda_{k})< \frac{\hat{\gamma}_{\nu}\eps^{2}}{L+4\eps}\right\}.
\end{equation}
Furthermore, let us define 
\begin{equation}\label{eq:K2}
\bK_{2}(\eps,x^{0},\nu,L):=\left\lceil\frac{(4\eps(\const-1/4)+f(x^{0})-f^{\ast})(L+4\eps)}{\hat{\gamma}_{\nu}\eps^{2}}\right\rceil
\end{equation}
and $\bK_{\max}(\eps,x^{0},\nu,L)=\min\{\bK_{1}(\eps,x^{0},\nu,L),\bK_{2}(\eps,x^{0},\nu,L)\}$. 

\begin{theorem}\label{th:main}
Let $(x^{k})_{k\geq 0}$ be generated by $\HBA(4\eps,L)$. Either the algorithm stops at $k_{\max}=\bK_{2}(\eps,x^{0},\nu,L)$, and reaches a point $x^{k_{\max}}\in\scrX^{\circ}$ satisfying 
\[
f(x^{k_{\max}})-f^{\ast}\leq\eps, 
\]
or else, it stops after $k_{\max}=\bK_{1}(\eps,x^{0},\nu,L)$ iterations, and we reach a point where $\lambda_{k_{\max}}\leq\eps$. 
\end{theorem}
\begin{Proof}
If $k_{\max}=\bK_{2}(\eps,x^{0},\nu,L) \leq \bK_{1}(\eps,x^{0},\nu,L)$, then for all $0\leq k\leq k_{\max}-1$ we have by definition $\tilde{\omega}_{\nu}(\lambda_{k})\geq \hat{\gamma}_{\nu}\frac{\eps^{2}}{L+4\eps}$. Hence, 
\begin{align*}
f^{\ast}&\leq f(x^{k_{\max}})\leq f(x^{0})-k_{\max}\hat{\gamma}_{\nu}\frac{\eps^{2}}{L+4\eps}+4\eps\const\leq f^{\ast}+\eps.
\end{align*}
For the second claim, suppose that $k_{\max}=\bK_{1}(\eps,x^{0},\nu,L)<\bK_{2}(\eps,x^{0},\nu,L)$, i.e. the algorithm stops before the objective function value is within $\eps$ of the global minimal value. Since $\lambda_{k}\to 0$, for all $\sigma>0$ there exists a $k$ such that $\lambda_{k}<\sigma$. Hence, let us fix a sufficiently small tolerance level  $\eps\in(0,1)$ so that the function $\tilde{\omega}_{\nu}$ is determined by terms including $\lambda_{k}^{2}$. Specifically, the following computations can be made for each generalized self-concordant parameter $\nu$: When $\nu\in(2,3)\cup(3,4)$ we see $\tilde{\omega}_{\nu}(\lambda_{k})= \frac{\tilde{\gamma}_{\nu}}{-\cb(L+4\eps)}\lambda^{2}_{k}$ for $k$ large enough, which smaller than $\frac{\tilde{\gamma}_{\nu}}{-\cb}\frac{\eps^{2}}{L+4\eps}$ exactly if $\lambda_{k}<\eps$. For $\nu\in\{3,4\}$ the same reasoning applies, proving the claim. 
\end{Proof}
\begin{remark}
Evaluating the stopping criterion $\bK_{1}$ appears to be expensive, since we have to keep track of the local norm of the search direction $\lambda_{k}$. However,  since $\tilde{\omega}_{\nu}$ is monotone, we can replace $\lambda_{k}$ with the more conservative figure $\sqrt{\tau_{h}}\beta_{k}$. Hence, if a bound on $\tau_{h}$ is available, we only have to monitor the evolution of the Euclidean length of the search direction.
\end{remark}

%%%%%%%%%%%%%%
While the above result is formulated in terms of convergence to stationary points of the potential function, our aim is actually to approximately solve the optimization problem \eqref{eq:P}. In order to connect these two conditions, we rely on our characterization of $\eps$-KKT points. Using the relation \eqref{eq:norm_equal}, we see 
\[
\norm{\nabla f(x^{k})-A^{\top}y^{k}}^{\ast}_{x^{k}}=\norm{\sqrt{H(x^{k})^{-1}}(\nabla f(x^{k})-A^{\top}y^{k})}_{2}\geq \frac{1}{\sqrt{\tau_{h}}}\norm{\nabla f(x^{k})-A^{\top}y^{k}}_{2}.
\]
Therefore, 
\begin{align*}
\norm{\nabla f(x^{k})-A^{\top}y^{k}}_{2}&\leq\sqrt{\tau_{h}}\norm{\nabla f(x^{k})-A^{\top}y^{k}}_{x^{k}}^{\ast}\\
&\leq \sqrt{\tau_{h}}\left(\norm{\nabla F_{\mu}(x^{k})-A^{\top}y^{k}}^{\ast}_{x^{k}}+\mu\norm{\nabla h(x)}^{\ast}_{x^{k}}\right)
\end{align*}
Recall that $\norm{\nabla F_{\mu}(x^{k})-A^{\top}y^{k}}^{\ast}_{x^{k}}=\lambda_{k}$. Furthermore, we know that $(x^{k})_{k\geq 0}\subseteq\scrS_{\mu}(x^{0})$, a compact set in $\scrX^{\circ}$. Since $h\in\bC^{3}(\scrC)$ and $\scrC$ contains no lines, the mapping, the norm $x\mapsto \norm{\cdot}^{\ast}_{x}$ is a continuous function on compact subsets of $\scrC$. Hence, the quantity 
\[
M_{\mu}(x^{0}):=\max_{x\in\scrS_{\mu}(x^{0})}\norm{\nabla h(x)}^{\ast}_{x},
\]
is well-defined and finite. In terms of this quantity we see that 
\begin{equation}\label{eq:boundfinal}
\norm{\nabla f(x^{k})-A^{\top}y^{k}}_{2}\leq \sqrt{\tau_{h}}\left(\lambda_{k}+\mu M_{\mu}(x^{0})\right),
\end{equation}
so that for $k\geq\bN(\eps/\sqrt{\tau_{h}},x^{0},\nu,L)$, and $\mu=\eps/\sqrt{\tau_{h}}$ we get 
\begin{align*}
\norm{\nabla f(x^{k})-A^{\top}y^{k}}_{2}=\bigoh(\eps).
\end{align*}
Combined with the inequalities \eqref{eq:sandwich1} we therefore conclude that $\chi(x^{k},y^{k})=\bigoh(\eps)$. i.e. we get and $\eps$-stationary point in the sense of Definition \ref{def:stationary}.  Note that $M_{\mu}(x^{0})$ is an algorithm independent constant, which can be computed before the method is started. Still it requires the solution of an optimization problem which can be fairly complicated in concrete instances, so it is definitely worthwhile searching for settings where this bound can be improved. Additionally, the complexity of the algorithm now explicitly depends on the eigenvalue bound $\tau_{h}$ of the barrier-generating kernel $h$, which means that if this number is big, the run time could become quite large.\footnote{However, both these remarks hold also for mirror descent type of methods, where the prox-function should be appropriately chosen since its properties affect the complexity bound}. Motivated by these observations, we next provide a refinement of this complexity result under the additional assumption that $h\in\scrF_{2,3}(\scrC)$ is a $\theta$-\emph{self-concordant barrier} in the sense of \eqref{eq:SC_barrier}. 

\begin{corollary}
\label{cor:scb_compl}
Let $\eps>0$ be a given tolerance level. If $h\in\scrF_{2,3}(\scrC)$ is a $\theta$-SCB, then running $\HBA(\eps/\sqrt{\theta},L)$ yields either an $2\sqrt{\tau_{h}}\eps$-stationary point, or an $\eps$ global minimum.
\end{corollary}
\begin{Proof}
For a $\theta$-SCB $h\in\scrF_{2,3}(\scrC)$, the complexity estimate in Theorem \ref{th:main} yields the following estimates: If 
\begin{align*}
k_{\max}=\bK_{1}(\eps,x^{0},3,L)=\min\left\{k\geq 0:\min\{\lambda_{k},\frac{\lambda_{k}^{2}}{L+\eps/\sqrt{\theta}}\}\leq \frac{\eps^{2}}{L+\eps/\sqrt{\theta}}\right\},
\end{align*}
then we know that the local norm of the gradient of the potential function is small, $\lambda_{k_{\max}}\leq\eps$. Then \eqref{eq:boundfinal} gives us 
\begin{align*}
\norm{\nabla f(x^{k_{\max}})-A^{\top}y^{k_{\max}}}_{2}\leq\sqrt{\tau_{h}}\left(\lambda_{k_{\max}}+\mu\sqrt{\theta}\right).
\end{align*}
Choosing $\mu=\eps/\sqrt{\theta}$, and using again the relation \eqref{eq:sandwich1}, the point $x^{k_{\max}}$ is seen to be a $2\sqrt{\tau_{h}}\eps$ stationary point, in the sense of Definition \ref{def:stationary}. If instead $k_{\max}=\bK_{2}(\eps,x^{0},3,L)$, we know we are $\eps$-close to the global minimum. 
\end{Proof}

\subsection{Analysis of $\AHBA(\mu)$}
The analysis of the adaptive version of our method follows similar lines as for the mother scheme $\HBA(\mu,L)$. The key innovation of the adaptive method is that it produces four recursive sequences $(x^{k})_{k\geq 0},(y^{k})_{k\geq 0},(\alpha_{k})_{k\geq 0},$ and $(L_{k})_{k\geq 0}$, where $\alpha_{k}=\alpha_{\mu}(x^{k},L_{k})$. We first show finite termination of the line search subroutine at each iteration, and establish a bound on the total number of function evaluations needed for its execution. The result is a generalization of the arguments in \cite{nesterov2006cubic,bogolubsky2016learning} for the case of relative smoothness in the non-convex case.

\begin{lemma}\label{lem:iteration}
Suppose that we run $\AHBA(\mu)$ for $N\geq 1$ rounds. Then, the total number of function evaluations $\scrE_{N}$, needed to satisfy \eqref{eq:line} in each of these $k=1,2,\ldots,N$ rounds, is at most 
\begin{equation}
\scrE_{N}\leq 2N+\log_{2}\left(\frac{2L}{L_{0}}\right).
\end{equation}
\end{lemma}
\begin{Proof}
Let $k=1,2,\ldots,N$ be an arbitrary iteration count. It is quite easy to see that the search cycle for $i_k$ is finite. Indeed since, by Definition \ref{def:L-adapt}, there exists such $L$ that for any $x,y\in\scrC$
\[
f(y) \leq f(x) + \inner{\nabla f(x),y-x}+LD_{h}(y,x),
\]
the search cycle for $i_k$ terminates no later than $L_{k+1}=2^{i_k-1}L_{k} \geq L$. At the same time, since $i_{k}\geq 0$ is the smallest integer for which \eqref{eq:line} holds, we have for $L_{k+1}/2=2^{i_{k}-2}L_{k}$ the inequality 
\begin{align*}
f(z^{k})>f(x^{k})+\inner{\nabla f(x^{k}),z^{k}-x^{k}}+2^{i_{k}-2}L_{k}D_{h}(z^{k},x^{k}).
\end{align*}
Hence, $2^{i_{k}-2}L_{k}<L$, or $L_{k+1}=2^{i_k-1}L_{k} \leq 2L$. Let us estimate the total number of function evaluations needed to ensure \eqref{eq:line}. On each iteration $k$, the number of function calls is $i_k+1=2+\log_2\frac{L_{k+1}}{L_{k}}$. Thus, the total number of function evaluations for $N$ rounds of execution of $\AHBA(\mu,L)$ is thus
\[
\scrE_{N}=\sum_{k=1}^N (i_k+1) = \sum_{k=1}^N \left(2+\log_2\frac{L_{k+1}}{L_{k}}\right) \leq 2N + \log_2\frac{2L}{L_0},
\] 
where we used the bound $L_k \leq 2L$. 
\end{Proof}
This shows that $\scrE_{N}=\bigoh(N)$, meaning that on average only a single function call is needed to satisfy the line search criterion \eqref{eq:line}. Thus, the performance of $\AHBA(\mu)$ is well described by the estimates for the overall iteration complexity of the method.

From the analysis of the base scheme $\HBA(\mu,L)$, we immediately deduce that the sequence $(x^{k})_{k\geq 0}$ generated by $\AHBA(\mu)$ satisfy the per-iteration descent 
\begin{align*}
F_{\mu}(x^{k+1})\leq F_{\mu}(x^{k})-\eta_{k}(x^{k},\alpha_{k})\equiv F_{\mu}(x^{k})-\Delta_{k},
\end{align*}
with the only difference that now the step size $\alpha_{k}$ is adaptively adjusted by evaluating the expression $\alpha_{\mu}(x^{k},2^{i_k-1}L_{k})$. From Lemma \ref{lem:omega}, we see that $\tilde{\omega}_{\nu}$ is a decreasing function of $L$. At the same time, $2^{i_k-1}L_{k}\leq 2L$ as it was shown above. This means that the adaptive versions of Lemma \ref{lem:complex_basic}, Theorem \ref{th:main}, and Corollary \ref{cor:scb_compl} are obtained by the change $L \to 2L$.
% Departing from here it is easy to verify that all the main results established for the base method $\HBA(\mu,L)$ hold as well. In particular, one can verify that the stopping time defined in \eqref{eq:K1} now becomes 
% $\bK_{1}(\eps,x^{0},\nu,4L)$, whereas, the stopping time \eqref{eq:K2} can stay as is. 
%Hence, the total gain thanks to the introduction of the adaptive choice of the $L$-smoothness constant is a speed-up of the method to the complexity estimate $\bK_{1}(\eps,x^{0},\nu,4L)$. This gain in efficiency promises significant gains over the base method. 
We see that the number of oracle calls increases for the adaptive version in comparison to non-adaptive. Nevertheless, the adaptive algorithm can use smaller values of $L$ and, hence, make longer steps, leading to faster convergence in practice.

%% file: SCAD.tex
%----------------------------------------------------------------------
%%% Example: SCAD
%----------------------------------------------------------------------
% !TEX root = ../Main.tex
%

We consider the non-convex statistical learning problem 
\begin{equation}\label{eq:folded}
\min_{\beta\in\R^{d}}\frac{1}{2}\norm{y-W\beta}^{2}_{2}+\sum_{i=1}^{d}p_{\zeta}(\abs{\beta_{i}}) 
\end{equation}
where $\ell(\beta):=\frac{1}{2}\norm{y-W\beta}^{2}_{2}$ is the quadratic data fitting term and $p_{\zeta}:\R_{+}\to\R_{+}$ is a \emph{folded concave penalty} \cite{LiuWanZha14,LohWain15,LiLiuYao16,LiLiuYaoYe17}, meaning that for given $a>2,\zeta>0$:
\begin{itemize}
\item[(i)] $t\mapsto p_{\zeta}(t)$ is non-decreasing and concave with $p_{\zeta}(0)=0$ and $p_{\zeta}(t)>0$ for $t>0$;
\item[(ii)] $t\mapsto p_{\zeta}(t)$ is differentiable on $[0,\infty)$;
\item[(iii)] $p'_{\zeta}(t)=0$ for all $t\geq a\zeta$ and $0\leq p'_{\zeta}(t)$ for any $t\geq 0$.
\end{itemize}
A specifc example would be \emph{smoothly clipped absolute deviation} (SCAD) penalty \cite{FanLi01} given by   
\begin{align*}
p_{\zeta}(t)=\left\{\begin{array}{ll} 
\zeta t & \text{if }0\leq t\leq \zeta,\\
\frac{1}{a-1}(-\frac{\zeta^{2}}{2}+a\zeta t-\frac{t^{2}}{2}) & \text{if } \zeta \leq t\leq a\zeta,\\
\frac{a+1}{2}\zeta^{2} & \text{if }t>a\zeta. 
\end{array}\right. 
\end{align*}
Note that the composite function $t\mapsto (p_{\zeta}\circ \abs{\cdot})(t)$ is continuous, but not differentiable at $t=0$. Hence, the objective function \eqref{eq:folded} is not smooth and non-convex. Doing some simple variable transformations, the regularized least-squares problem \eqref{eq:folded} can be put into an optimization problem fitting the structure of this paper. Let us introduce new variables  $\beta^{+}_{i}:=\max\{\beta_{i},0\}$ and $\beta_{i}^{-}:=\max\{-\beta_{i},0\}$, so that $\beta_{i}^{+}+\beta_{i}^{-}=\abs{\beta_{i}}$. We additionally allow the inclusion of a-priori upper bounds on the parameter vector. This gives rise to a box-constrained reformulation of \eqref{eq:folded} of the form 
\begin{align*}
 & \min_{\beta^{+}\in\R^{d},\beta^{-}\in\R^{d}}\ell(\beta^{+}-\beta^{-})+\sum_{i=1}^{d}p_{\zeta}(\beta_{i}^{+}+\beta_{i}^{-}),\\
\text{s.t. } & 0\leq \beta_{i}^{-}\leq u_{i},\; 0\leq \beta_{i}^{+}\leq u_{i}\quad 1\leq i\leq d. 
\end{align*}
To bring this problem into a formulation fitting this paper, we first relabel the pair $(\beta^{-},\beta^{+})\in\R^{d}\times\R^{d}$ into one long vector $x:=(x_{1},\ldots,x_{d},x_{d+1},\ldots,x_{2d})$, where the first $d$ variables correspond to the positive part and the remaining $d$ variables represent the negative part. Call $n:=2d$ we define the data fitting term to be $f_{0}(x):=\ell(Bx)$, where $B:\R^{n}\to\R^{d}$ is the linear operator $(Bx)_{i}:=x_{d+i}-x_{i}$ for all $i\in\{1,2,\ldots,d\}$. The regularizer can be written as $f_{1}(Dx):=\sum_{i=1}^{d}p_{\zeta}(x_{i}+x_{d+i})$, corresponding the the composition of the function $\R^{d}\ni y\mapsto f_{1}(y)=\sum_{i=1}^{d}p_{\zeta}(y_{i})$ with the linear operator $D:\Rn\to\R^{d}$ given by $(Dx)_{i}=x_{i}+x_{d+i}$ for all $i\in\{1,2,\ldots,d\}$. Define $\scrX=\bar{\scrC}:=\prod_{i=1}^{n}[0,u_{i}]$, so that our non-convex minimization problem reads as 
\begin{equation}\label{eq:Regularized}
\min_{x\in\scrX}\{f(x):=f_{0}(x)+f_{1}(Dx)\},
\end{equation}
where $f_{0}(x):=\frac{1}{2}x^{\top}Qx+x^{\top}q$ is a convex quadratic function with Hessian $Q:=B^{\top}W^{\top}WB$ and $q^{\top}:=-B^{\top}W^{\top}y$. Note that $\bar{\scrC}$ admits a simple self-concordant function (e.g. the Burg entropy as described in Example \ref{ex:bgf}), but is not prox-friendly (see Remark \ref{rem:box}).

\begin{figure}
    \centering
   \subfigure[Fitted Values over the test sample]{
        \includegraphics[width=0.4\textwidth]{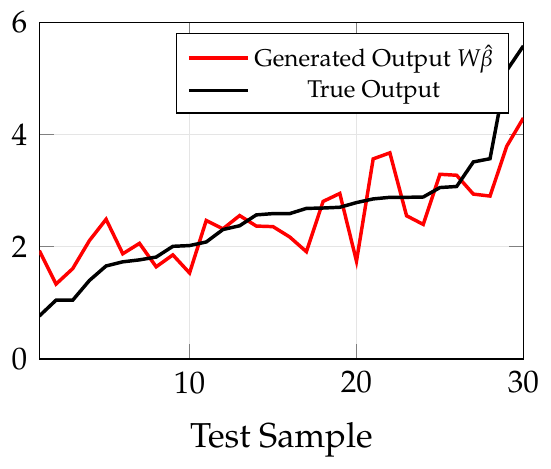}
      % \caption{}
        }
     \subfigure[Test Error of $\HBA(\mu,L)$]{
        \includegraphics[width=0.4\textwidth]{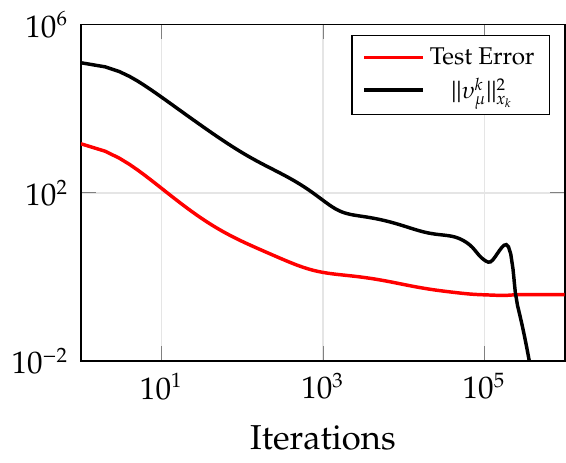}
     %   \caption{}
     }
    \caption{Fitted values and Test error of $\HBA(\mu,0)$ for the Prostate Cancer data with SCAD regularization. As parameters for the SCAD regularizer we have used $\zeta = 0.01$, $a=10$. For the barrier parameter the value $\mu=10^{-3}$ has been chosen.}
    \label{fig:StatLearn}
\end{figure}

The quadratic loss function $f_{0}(x)$ is convex and has a Lipschitz continuous gradient with Lipschitz constant $\rho:=\abs{Q}$. Rescaling the data appropriately, we can assume without loss of generality that $\rho\geq 1$. Hence, for the data fidelity part, a standard Lipschitz-descent lemma \cite{NesConvex} applies and gives 
\begin{equation}\label{eq:Lip1}
f_{0}(y)\leq f_{0}(x)+\inner{\nabla f_{0}(x),y-x}+\frac{\rho}{2}\norm{y-x}^{2}_{2}. 
\end{equation}
For $\theta\in\R^{d}_{+}$, the regularizing term reads as $f_{1}(\theta)=\sum_{i=1}^{d}p_{\zeta}(\theta_{i})$, and each summand in this expression is a concave and differentiable function on $(0,\infty)$. Hence, for all $s,t>0$, we have
\begin{align*}
p_{\zeta}(s)\leq p_{\zeta}(t)+p'_{\zeta}(t)(s-t).
\end{align*}
For any two vectors $\theta^{(1)},\theta^{(2)}\in\R^{d}_{++}$ this implies that
\begin{align*}
f_{1}(\theta^{(2)})\leq f_{1}(\theta^{(1)})+\inner{\nabla f_{1}(\theta^{(1)}),\theta^{(2)}-\theta^{(1)}}.
\end{align*}
Evaluating this expression at the vectors $\theta^{(1)}=Dx$ and $\theta^{(2)}=Dy$, we obtain 
\begin{equation}\label{eq:Lip2}
f_{1}(Dy)\leq f_{1}(Dx)+\inner{D^{\top}\nabla f_{1}(Dx),y-x}.
\end{equation}
Adding \eqref{eq:Lip1} with \eqref{eq:Lip2}, we see that 
\begin{equation}\label{eq:descentfinal}
f(y)\leq f(x)+\inner{\nabla f(x),y-x}+\frac{\rho}{2}\norm{y-x}^{2}_{2}\qquad\forall y,x\in\scrX^{\circ}.
\end{equation}
For the rest of the analysis we assume that $u_{i}=\infty$, so that no external upper bounds on the parameter vectors are imposed. Thus, $\scrX=\bar{\scrC}=\Rn_{+}$, and the natural barrier-generating kernel for this set is the Burg entropy $h(x)=-\sum_{i=1}^{n}\ln(x_{i})$, inducing the Riemannian metric $H(x)=\diag\{x_{1}^{-2},\ldots,x_{n}^{-2}\}$, and Bregman divergence 
\[
D_{h}(y,x)=h(y)-h(x)-\inner{\nabla h(x),y-x}=\sum_{i=1}^{n}\ln\left(\frac{y_{i}}{x_{i}}\right)+\sum_{i=1}^{n}\frac{y_{i}}{x_{i}}-n. 
\]
In terms of the potential function $F_{\mu}(x)=f(x)+\mu h(x)$, the combined descent inequality \eqref{eq:descentfinal} reads as 
\[
F_{\mu}(y)\leq F_{\mu}(x)+\inner{\nabla F_{\mu}(x),y-x}+\frac{\rho}{2}\norm{y-x}_2^2+\mu D_{h}(y,x).
\]
Defining the regularized Burg entropy $\tilde{h}_{\mu,\rho}:=-\sum_{i=1}^{n}\ln(x_{i})+\frac{\rho}{2\mu}\norm{x}^{2}_{2}$, we can write the descent inequality for the potential function in more concise terms as 
\[
F_{\mu}(y)\leq F_{\mu}(x)+\inner{\nabla F_{\mu}(x),y-x}+\mu D_{\tilde{h}_{\mu,\rho}}(y,x).
\]
Note that $\tilde{h}_{\mu,\rho}\in\scrF_{2,3}(\R^{n}_{++})$. This shows that the regularized statistical learning problem can be solved with $\HBA(\mu,0)$. We apply the model to the Prostate Cancer data set available at \url{https://web.stanford.edu/~hastie/ElemStatLearn/data.html}. This data set consists of a total of 97 samples with 8 dimensions each, from which 67 are used to train the model and 30 are used for validation. Thus, in this case we have a matrix $W\in\mathbb{R}^{67\times8}$ and $y \in \mathbb{R}^{67}$. Moreover, we have used the following set of parameter values $\zeta = 0.01$, $a=10$, and $\mu = 1\cdot10^{-3}$. Once a model $\hat \beta$ is found, such value is used to predict a output $A_{\text{test}\hat \beta}$ for the test database $A_{\text{test}}$. Figure~\ref{fig:StatLearn}(a) shows the true output of the test database for each of the $30$ data points colored in black, and the predicted output for the same points generated by the output of the Adaptive HBA algorithm. Moreover, Figure~\ref{fig:StatLearn}(b) shows the gradient norm value versus the number of iterations of the algorithm, and the test error in color red. $\AHBA(\mu)$ reaches a test error of $0.363$.  This value improves upon the $0.4194$ test error reported in \cite{BiaCheYe15}, and the $0.479$ test error reported in \cite[Table 3.3]{FriHasTib01}.

%% file: Lp.tex
%----------------------------------------------------------------------
%%% Example: Lp
%----------------------------------------------------------------------
% !TEX root = ../Main.tex
%

Consider the optimization problem 
\begin{equation}\label{eq:BNP}
\begin{array}{ll}
\min &  f(x)=\sum_{i=1}^{n}x_{i}^{p} \\
\textnormal{subject to} & x\in\scrX=\R_{+}^{n}\cap\scrA
\end{array}
\end{equation}
where the problem inputs consist of $A\in\R^{m\times n}$, $b\in\R^{m}$ and $p\in(0,1]$. Sparse signal or solution reconstruction by solving problem \eqref{eq:BNP}, especially for the case where $p\in(0,1)$, has recently received considerable attention; see e.g. \cite{BruDonEla09}. In signal reconstruction, one typically has linear measurements $b=Ax$, where $x$ is a sparse signal, i.e. the sparsest or smallest support cardinality solution of the linear system. This sparse signal is recovered by solving the inverse problem \eqref{eq:BNP} with the non-smooth, non-convex objective function $\norm{x}_{0}=\abs{\{i\in\{1,2,\ldots,n\}\vert x_{i}>0\}}$. The $L_{0}$-norm optimization problem is shown to be NP-hard. When $p=1$, the problem is reduced to a linear program, and hence it can be solved in polynomial time. If $p>1$, the problem \eqref{eq:BNP} becomes a convex optimization optimization problem, and thus is also efficiently solvable with fast interior point methods. Only recently, the challenging case where $p\in(0,1)$ has been thoroughly investigated in \cite{GeJiaYe11,BiaCheYe15}. We aim to solve this NP-hard problem with $\AHBA(\mu)$. Given the geometry, it is natural to look at the barrier-generating kernel $h(x)=-\sum_{j=1}^{n}\ln(x_{j})$, so that the potential function $F_{\mu}$ becomes $F_{\mu}(x)=\norm{x}^{p}_{p}-\mu\sum_{i=1}^{n}\ln(x_{i})$. Note that the objective function $f(x)=\norm{x}^{p}_{p}$ is twice continuously differentiable on $\scrX^{\circ}$ and concave. Hence, L-smoothness holds for any $L>0$.

\begin{figure}
    \centering
   \subfigure[Fitted Values over the test sample]{
        \includegraphics[width=0.4\textwidth]{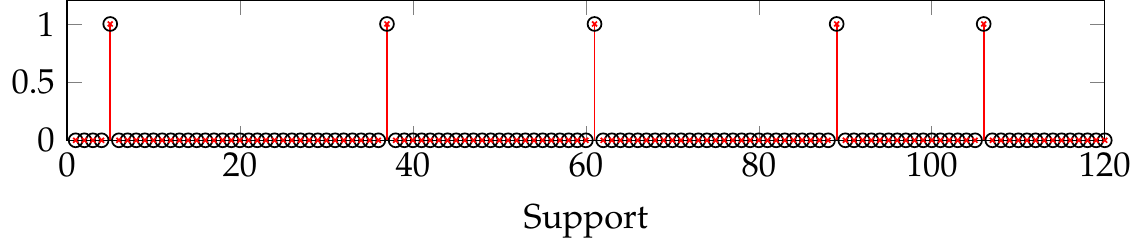}
      % \caption{}
        }
     \subfigure[Recovered Signal]{
        \includegraphics[width=0.4\textwidth]{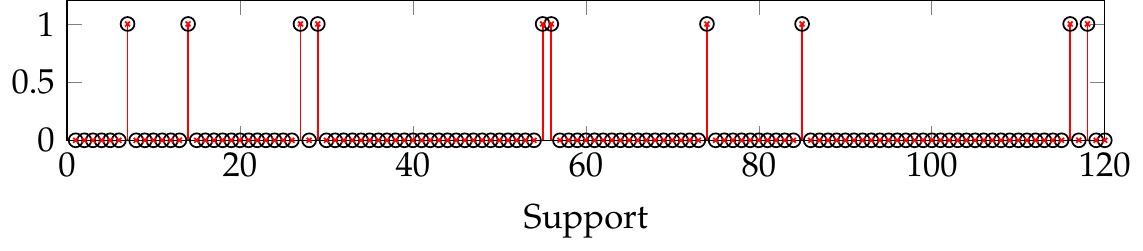}
     %   \caption{}
     }
         \subfigure[Success Rate of $\HBA(\mu,L)$]{
        \includegraphics[width=0.4\textwidth]{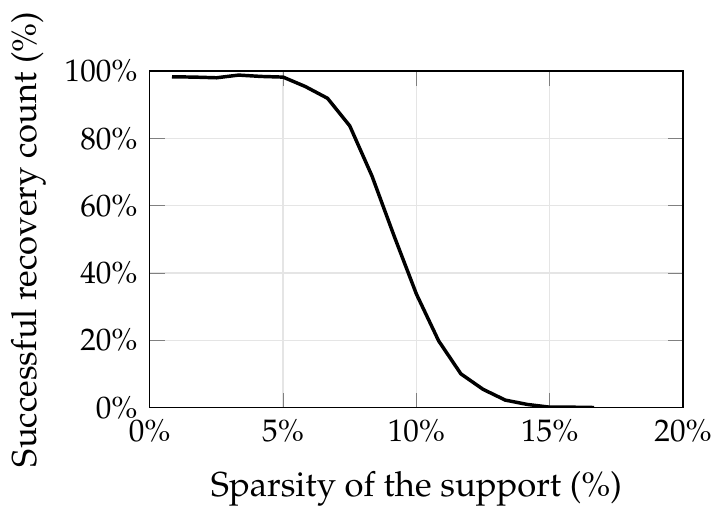}
     %   \caption{}
     }
    \caption{Recovery pattern and recovery rates for the $L^{p}$ minimization problem solved by $\AHBA(\mu)$. As parameter values we have chosen $p=0.5$ and $\mu=1$. }
    \label{fig:LpEstimation}
\end{figure}

To test the performance of our method, we have set up numerical experiments and recorded the recovery rate of the true underlying signal for each level of sparsity. Specifically, we generate a binary signal of length $120$, denoted as $\hat x$, and various sparsity patterns. The excellent recovery properties of $AHBA(\mu)$ with $5$ non-zero entries is displayed in Figure~\ref{fig:LpEstimation}(a) and with $10$ non-zero entries in~Figure\ref{fig:LpEstimation}(b), in which the original signal is marked as black circles $\circ$ and the recovered one is mark as red crosses {\color{red}$\times$}. In each case, we generated an observation matrix as an orthogonal positive sensing matrix $A$, and a set of $30$ observations. Moreover we have used as parameter values $p=0.5$ and $\mu=1$. Figure \ref{fig:LpEstimation}(c) reveals the general pattern of the recovery rates of the true signal.

%% file: conclusion.tex
%%% Conclusion
%----------------------------------------------------------------------
% !TEX root = ./Main.tex
%%%

In this paper we introduced and studied a new class of interior point methods based on the Hessian-barrier technique originally developed in \cite{HBA-linear}. Using the theory of generalized self-concordant functions we can significantly extend the applicability of this method to cover general non-convex optimization problems on polyhedral domains with a set constraint admitting a generalized self-concordance set-up. Theoretical convergence and complexity results are proven, showing that the method achieves the optimal iteration complexity $O(\eps^{-2})$. We have tested the method empirically and verified that the method performs also well in practice. There are many important directions for future investigations to be made. First, it is very important to relax the present algorithmic scheme to allow for inexact computations and to derive a path-following approach allowing the barrier parameter $\mu$ to vary over the run time of the algorithm. Allowing for numerical and random noise is of relevance when HBA methods are to be designed in distributed optimization settings \cite{LinZha15,GasLeeNedUri18}, something we plan to do in the future, and the path-following approach might allow us to strengthen the convergence properties of the algorithm.

The most costly step of $\HBA(\mu,L)$, and its adaptive version $\AHBA(\mu)$, is the solution of the linear system \eqref{eq:finder}. If the function $f$ appears as a finite sum, a very important direction for future research is to either use preconditioning techniques or randomization and sketching ideas, to speed up the computation. Second, the method should be also a competitive first-order scheme for stochastic optimization. In fact, self-concordant functions have been already successfully used in online learning \cite{NarRak10} and random sampling \cite{Nar16}, and the class of generalized self-concordant functions may provide interesting extensions of these seminal contributions. Finally, it will be important to identify acceleration strategies for the basic HBA template to make it even more attractive for large-scale application in engineering and machine learning. This is another challenging line of research we plan to pursue in the near future.  

%% file: Appendix_Step.tex
%%% Appendix_Step
%----------------------------------------------------------------------
% !TEX root = ./Main.tex
%

For each parameter $\nu\in(2,4]$ we derive the corresponding optimal step-size policy via a simple optimization argument. This will prove the claimed optimality of the policy. 
%%%%%%%%%%%%%%%%%%%%%%%%%%%%
\subsection{The case $\nu\in(2,3)$}
Using the definition 
\[
\omega_{\nu}(t)=\left(\frac{\nu-2}{4-\nu}\right)\frac{1}{t}\left[\frac{\nu-2}{2t(3-\nu)}\left((1-t)^{\frac{2(3-\nu)}{2-\nu}}-1\right)-1\right],
\]
we obtain 
\begin{align*}
\eta_{\mu}(x,t)=&t\lambda^{2}_{\mu}(x)-t^{2}\lambda^{2}_{\mu}(x)(L+\mu)\frac{(\nu-2)^{2}}{2t^{2}\delta_{\mu}(x)^{2}(3-\nu)(4-\nu)}\left((1-t\delta_{\mu}(x))^{\frac{2(3-\nu)}{2-\nu}}-1\right)\\
&+t\frac{\lambda^{2}_{\mu}(x)}{\delta_{\mu}(x)}(L+\mu)\frac{\nu-2}{4-\nu}\\
&=t\left(\lambda^{2}_{\mu}(x)+\lambda^{2}_{\mu}(x)\frac{L+\mu}{\delta_{\mu}(x)}\frac{\nu-2}{4-\nu}\right)\\
&-\left(\frac{\lambda_{\mu}(x)}{\delta_{\mu}(x)}\right)^{2}(L+\mu)\frac{(\nu-2)^{2}}{2(3-\nu)(4-\nu)}\left((1-t\delta_{\mu}(x))^{\frac{2(3-\nu)}{2-\nu}}-1\right).
\end{align*}

For $t\delta_{\mu}(x)\in(0,1)$, this gives 
\begin{align*}
&\frac{\partial}{\partial t}\eta_{\mu}(x,t)=\lambda^{2}_{\mu}(x)\left(1+\frac{L+\mu}{\delta_{\mu}(x)}\frac{\nu-2}{4-\nu}\right)-\frac{\lambda^{2}_{\mu}(x)}{\delta_{\mu}(x)}\frac{(L+\mu)(\nu-2)}{4-\nu}(1-t\delta_{\mu}(x))^{-\frac{4-\nu}{\nu-2}},\\
&\frac{\partial^{2}}{\partial t\partial t}\eta_{\mu}(x,t)=-\lambda^{2}_{\mu}(x)(L+\mu)(1-t\delta_{\mu}(x))^{\frac{2}{2-\nu}}<0. 
\end{align*}
Solving the stationarity condition $\frac{\partial}{\partial t}\vert_{t=\alpha_{\mu}(x,L)}\eta_{\mu}(x,t)=0$, gives 
\begin{equation}
\alpha_{\mu}(x,L)\delta_{\mu}(x)=1-\left(1+\frac{\delta_{\mu}(x)}{L+\mu}\frac{4-\nu}{\nu-2}\right)^{-\frac{\nu-2}{4-\nu}}. 
\end{equation}
Since $\frac{\nu-2}{4-\nu}\in(0,1)$ for $\nu\in(2,3)$, the Bernoulli inequality gives 
\begin{align*}
\left(1+\frac{\delta_{\mu}(x)}{L+\mu}\frac{4-\nu}{\nu-2}\right)^{\frac{\nu-2}{4-\nu}}\leq 1+\frac{\delta_{\mu}(x)}{L+\mu},
\end{align*}
so that, for $L\geq 0$,
\begin{align*}
1-\left(1+\frac{\delta_{\mu}(x)}{L+\mu}\frac{4-\nu}{\nu-2}\right)^{-\frac{\nu-2}{4-\nu}}\leq \frac{\delta_{\mu}(x)}{\delta_{\mu}(x)+(L+\mu)}<1. 
\end{align*}
Hence, setting 
\begin{equation}\label{eq:step_mid}
\alpha_{\mu}(x,L)=\frac{1}{\delta_{\mu}(x)}\left[1-\left(1+\frac{\delta_{\mu}(x)}{L+\mu}\frac{4-\nu}{\nu-2}\right)^{-\frac{\nu-2}{4-\nu}}\right]
\end{equation}
gives $\eta_{\mu}(x,\alpha_{\mu}(x,L))>\eta_{\mu}(x,0)=0$ and $\metric_{\nu}(x,T_{\mu}(x,\alpha_{\mu}(x,L))<1$.  
%%%%%%%%%%%%%%%%%%%%%%%%%%%
\subsection{The case $\nu=3$.}
We have $\omega_{3}(t)=\frac{-1}{t^{2}}\left(t+\ln(1-t)\right)$ for $t\in(-\infty,1)$, and $\metric_{3}(x,y)=\frac{M}{2}\norm{y-x}_{x}$. Hence, $\delta_{\mu}(x)=\frac{M}{2}\lambda_{\mu}(x)$, and 
\begin{align*}
\eta_{\mu}(x,t)=t\lambda^{2}_{\mu}(x)+\left(\frac{\lambda_{\mu}(x)}{\delta_{\mu}(x)}\right)^{2}(L+\mu)\left[t\delta_{\mu}(x)+\ln(1-t\delta_{\mu}(x))\right]. 
\end{align*}
Therefore,
\begin{align*}
&\frac{\partial}{\partial t}\eta_{\mu}(x,t)=\lambda_{\mu}(x)^{2}+\left(\frac{\lambda_{\mu}(x)}{\delta_{\mu}(x)}\right)^{2}(L+\mu)\left[\delta_{\mu}(x)-\frac{\delta_{\mu}(x)}{1-t\delta_{\mu}(x)}\right]\\
&\frac{\partial^{2}}{\partial t\partial t}\eta_{\mu}(x,t)=-\frac{\lambda_{\mu}(x)^{2}}{\delta_{\mu}(x)}(L+\mu)(1-t\delta_{\mu}(x))^{-2}<0
\end{align*}
Solving for the stationary condition $\frac{\partial}{\partial t}\vert_{t=\alpha_{\mu}(x)}\eta_{\mu}(x,t)=0$, we get
\begin{align*}
1-\alpha_{\mu}(x)\delta_{\mu}(x)=t(L+\mu).%\iff \alpha_{\mu}(x)=\frac{1}{\delta_{\mu}(x)+L+\mu}. 
\end{align*}
Hence, setting 
\begin{equation}
\alpha_{\mu}(x)=\frac{1}{\delta_{\mu}(x)+L+\mu},
\end{equation}
we observe that $\eta_{\mu}\left(x,\alpha_{\mu}(x,L)\right)>\eta_{\mu}(x,0)=0$, and $\alpha_{\mu}(x,L)\delta_{\mu}(x)<1$. Therefore, $\metric_{3}\left(x,T_{\mu}(x,\alpha_{\mu}(x,L))\right)<1$.
%%%%%%%%%%%%%%%%%%%%%%%%%%%%%%%%%%%%%%%
\subsection{The case $\nu=4$.}
We have $\omega_{4}(t)=\frac{(1-t)\ln(1-t)+t}{t^{2}}$, and $\metric_{4}(x,y)=M\norm{y-x}^{-1}_{2}\norm{y-x}_{x}^{2}$. Hence, $\delta_{\mu}(x)=M\frac{\lambda^{2}_{\mu}(x)}{\beta_{\mu}(x)}$, and 
\begin{align*}
\eta_{\mu}(x,t)=t\lambda_{\mu}(x)^{2}-\left(\frac{\lambda_{\mu}(x)}{\delta_{\mu}(x)}\right)^{2}(L+\mu)[t\delta_{\mu}(x)+(1-t\delta_{\mu}(x))\ln(1-t\delta_{\mu}(x))].
\end{align*}
Therefore, 
\begin{align*}
&\frac{\partial}{\partial t}\eta_{\mu}(x,t)=\lambda_{\mu}(x)^{2}-\left(\frac{\lambda_{\mu}(x)}{\delta_{\mu}(x)}\right)^{2}(L+\mu)\ln(1-t\delta_{\mu}(x)),\\
&\frac{\partial^{2}}{\partial t\partial t}\eta_{\mu}(x,t)=-\frac{\lambda^{2}_{\mu}(x)(L+\mu)}{1-t\delta_{\mu}(x)}<0. 
\end{align*}
Solving for stationarity $\frac{\partial}{\partial t}\vert_{t=\alpha_{\mu}(x,L)}\eta_{\mu}(x,t)=0$, gives 
\begin{align*}
\frac{-\delta_{\mu}(x)}{L+\mu}=\ln\left(1-\alpha_{\mu}(x,L)\delta_{\mu}(x,L)\right),
\end{align*}
so that 
\begin{equation}
\alpha_{\mu}(x,L)=\frac{1}{\delta_{\mu}(x,L)}\left[1-\exp\left(-\frac{\delta_{\mu}(x)}{L+\mu}\right)\right].
\end{equation}
It follows $\metric_{4}(x,T_{\mu}(x,\alpha_{\mu}(x,L)))=1-\exp(-\delta_{\mu}(x)/(L+\mu))\in(0,1)$. 

%%%%%%%%%%%%%%%%%%%%%%%%%%%%%%%%%%%%%%%%%%%%%%%%%%%%%%%%%%%%%
\subsection{The case $\nu\in(3,4)$.} 
The basic computations for this range can be copied from the case $\nu\in(2,3)$. Doing so, we immediately arrive at the step size policy
\begin{equation}
\alpha_{\mu}(x,L)\delta_{\mu}(x)=1-\left(1+\frac{\delta_{\mu}(x)}{L+\mu}\frac{4-\nu}{\nu-2}\right)^{-\frac{\nu-2}{4-\nu}}. 
\end{equation}
From here, we can continue all the computations as for the case $\nu\in (2,3)$ to conclude that the step size $\alpha_{\mu}(x,L)$ is given by \eqref{eq:step_mid}. Note that $\frac{2-\nu}{4-\nu}<0$, so that $(1+\frac{\delta_{\mu}(x)}{L+\mu}\frac{4-\nu}{\nu-2})^{\frac{2-\nu}{4-\nu}}\in(0,1)$, and therefore $\alpha_{\nu}(x)\delta_{\mu}(x)\in(0,1)$. All other conclusions derived for $\nu\in(2,3)$ apply to the present setting as well. 
%%%%%%%%%%%%%%%%%%%%%%%%%%%%%%%%%%%%%%%%%%%%%%%

%% file: appendix_stationary.tex
%%% Proof Stationary 
%--------------------------------------------------------------------
% !TEX root = ./Main.tex
%

We denote by $\Delta_{k}\equiv \eta_{\mu}(x^{k},\alpha_{k})$, where $x^{k}$ is the iterate of $\HBA(\mu,L)$, and $\alpha_{k}\equiv \alpha_{\mu}(x^{k},L)$ is the associated step size. Similarly, we define the sequence $\lambda_{k},\beta_{k}$ and $\delta_{k}$ as in \eqref{eq:lambda_beta_delta}. %Recall that 
%\[
%\delta_{k}=M\left(\frac{\nu}{2}-1\right)\lambda_{k}^{\nu-2}\beta_{k}^{3-\nu}\qquad\forall k\geq 0.
%\]

\subsection{The case $\nu\in(2,3)$}
An explicit calculation shows that 
\begin{align*}
\Delta_{k}&=\frac{\lambda^{2}_{k}}{\delta_{k}}\left[1-\frac{4-\nu}{2(3-\nu)}\left(1+\frac{\delta_{k}}{L+\mu}\frac{4-\nu}{\nu-2}\right)^{(2-\nu)/(4-\nu)}\right]\\
&+\left(\frac{\lambda_{k}}{\delta_{k}}\right)^{2}\frac{(\nu-2)(L+\mu)}{2(3-\nu)}\left[1-\left(1+\frac{\delta_{k}}{L+\mu}\frac{4-\nu}{\nu-2}\right)^{(2-\nu)/(4-\nu)}\right]. 
\end{align*}
To make the analysis of this expression more convenient, we introduce the quantities  
\begin{align*}
t_k &:= 1-\frac{1}{\cb}\frac{\delta_{k}}{L+\mu} \in (1,+\infty),\text{ and }\\
\ca&:=\frac{4-\nu}{2(3-\nu)}\in(1,+\infty),\; \cb:=\frac{2-\nu}{4-\nu}\in(-1,0). 
\end{align*}
Then $\frac{(\nu-2)(L+\mu)}{\delta_k}=\frac{(4-\nu)}{t_k-1}$ and 
\begin{align*}
\Delta_{k}&=\frac{\lambda^{2}_{k}}{\delta_{k}} \left( 1-\ca t_k^\cb + \frac{\ca}{t_k-1}(1-t_k^\cb)\right) \\
&= \frac{\lambda^{2}_{k}}{\delta_{k}} \left( 1+\frac{\ca}{t_k-1}-\ca t_k^\cb\left(1+\frac{1}{t_k-1} \right) \right) \\
&= \frac{\lambda^{2}_{k}}{\delta_{k}} \left(1+\frac{\ca}{t_k-1}-\frac{\ca t_k^{\cb+1}}{t_k-1}\right). 
\end{align*} 
Let us define a function $\gamma(t)$ such that $\Delta_{k}=\frac{\lambda^{2}_{k}}{\delta_{k}}\gamma(t_k)$. 
Our next goal is to show that, for $t \in [2,+\infty)$, $\gamma(t)$ is below bounded by some positive constant and, for $t \in (1,2]$, $\gamma(t)$ is below bounded by some positive constant multiplied by $t-1$.

\textit{1. $t \in [2,+\infty)$.} We will show that $\gamma'(t) \geq 0$, whence $\gamma(t) \geq \gamma(2)$. Thus, we need to show that
\begin{align*}
0 \leq & %\leq \frac{d}{dt} \left( 1+\frac{\ca}{t-1}-\frac{\ca t^{\cb+1}}{t-1}\right) 
\gamma'(t)= -\frac{\ca }{(t-1)^2} \underbrace{\left(1-(\cb+1)t^\cb+\cb t^{\cb+1}\right)}_{=:\psi(t)}. 
\end{align*}
Since $\ca>1$, to show that $\gamma'(t) \geq 0$ it is enough to show that $\psi(t) \leq 0$. Since $\cb \in(-1,0)$ and $t\geq 2$,
\[
\psi'(t) = \cb(\cb+1)t^\cb- \cb(\cb+1) t^{\cb-1} = \cb(\cb+1)t^{\cb-1}(t-1) \leq 0.
\]
Whence, $\psi(t) \leq \psi(2)$ for all $t \in [2,+\infty)$. It remains to show that $\psi(2) \leq 0$. Let us consider $\psi(2) =\varphi(\cb):= 1-(\cb+1)2^\cb+\cb 2^{\cb+1} = 1+\cb 2^\cb-2^\cb$ as a function of $\cb\in(-1,0)$. Clearly, $\varphi(-1)=\varphi(0)=0$, and it is easy to check via the intermediate value theorem that $\varphi(b)<0$ for all $b\in(-1,0)$. We conclude that for $t\geq 2$ we get $\Delta_{k}\geq \frac{\lambda^{2}_{k}}{\delta_{k}}\gamma(2)$.

\textit{2. $t \in (1,2]$.} We will show that $ \frac{d}{dt} \left(\gamma(t)/(t-1)\right) \leq 0$, whence $\gamma(t) \geq (t-1)\gamma(2)$. Thus, we need to show that
\begin{align*}
0 & \geq \frac{d}{dt} \left( \frac{1}{t-1}+\frac{\ca}{(t-1)^2}-\frac{\ca t^{\cb+1}}{(t-1)^2}\right) \\
%&= -\frac{1}{(t-1)^2}-\frac{2\ca}{(t-1)^3}-\frac{\ca(\cb+1)t^\cb}{(t-1)^2}+\frac{2\ca t^{\cb+1}}{(t-1)^3} \\
%& = \frac{1}{(t-1)^3} \left( -t + 1 -2 \ca - \ca(\cb+1)t^{\cb+1} + \ca(\cb+1)t^\cb + 2\ca t^{\cb+1} \right) \\
& = \frac{1}{(t-1)^3} \left( -t + 1 -2 \ca + \ca(\cb+1)t^{\cb} - \ca(\cb-1)t^{\cb+1} \right)\equiv \frac{1}{(t-1)^{3}}\psi(t). 
\end{align*}
Therefore, our next step is to show that $\psi(t) \leq 0$. We have
\begin{align*}
\psi'(t) & =-1 + \ca (\cb+1) \cb t^{\cb-1}-\ca(\cb-1)(\cb+1) t^{\cb},\\
\psi''(t) & =\ca\cb (\cb+1) (\cb-1) t^{\cb-2}-\ca(\cb-1)\cb(\cb+1) t^{\cb-1} \\
& = \ca \cb(\cb+1) (\cb-1) t^{\cb-2} (1-t).
\end{align*}
By definition, $\ca(\cb+1) = 1$. Hence, since $t >1$ and $\cb\in(-1,0)$, we observe that $\psi''(t) \leq 0$. Thus, $\psi'(t) \leq \psi'(1) = 0$, and consequently, $\psi(t) \leq \psi(1) = 0$, for all $t\in(1,2]$. This proves the claim $\gamma(t)/(t-1) \geq \gamma(2)$ for $t\in(1,2]$.

Combining both cases, we obtain that $\gamma(t) \geq \min\{\gamma(2), (t-1) \gamma(2)\}$, where $\gamma(2)=1-\ca+\ca 2^{1/\ca}$, using the fact that $\cb+1=1/\ca$. Unraveling this expression by using the definition of the constant $\ca$, we see that $\gamma(2)$ depends only on the self-concordance parameter $\nu\in(2,3)$. In light of this, let us introduce the constant 
\begin{equation}\label{eq:tildegamma}
\tilde{\gamma}_{\nu}:=1+\frac{4-\nu}{2(3-\nu)}\left(1-2^{2(3-\nu)/(4-\nu)}\right).
\end{equation}
Observe that $\tilde{\gamma}_{2}=0$ and, by a simple application of l'H\^{o}pital's rule, $\lim_{\nu\uparrow 3}\hat{\gamma}_{\nu}=1-\log(2)\in(0,1)$.
Hence $\gamma(2)\equiv\tilde{\gamma}_{\nu}\in(0,1)$ for all $\nu\in(2,3)$. We conclude, 
\begin{equation*}
\Delta_k\geq \frac{\tilde{\gamma}_{\nu} \lambda^{2}_{k}}{\delta_{k}} \min\left\{1,\frac{-1}{\cb} \frac{\delta_{k}}{L+\mu}\right\} = \tilde{\gamma}_{\nu} \min\left\{\frac{\lambda^{2}_{k}}{\delta_{k}}, \frac{\lambda^{2}_{k}}{L+\mu}\frac{-1}{\cb}\right\}.
\end{equation*}
Since $\lambda_{k}\geq \sqrt{\sigma_{h}}\beta_{k}$, $ \delta_{k}=M(\frac{\nu}{2}-1)\lambda_{k}^{\nu-2}\beta_{k}^{3-\nu}$, 
the following lower and upper bounds can be established for $\nu\in(2,3)$:
\begin{align*}
M\left(\frac{\nu}{2}-1\right)\sigma_{h}^{\frac{\nu-2}{2}}\alpha_{k}\beta_{k}\leq\alpha_{k}\delta_{k}\leq M\left(\frac{\nu}{2}-1\right)\sigma_{h}^{-\frac{3-\nu}{2}}\alpha_{k}\lambda_{k}.
\end{align*}
This estimate implies first that $\delta_{k}\leq M(\nu/2-1)\sigma_{h}^{-(3-\nu)/2}\lambda_{k} $, and second
\begin{align*}
\frac{\lambda_{k}^{2}}{\delta_{k}}\geq\frac{2\lambda_{k}}{M(\nu-2)}\sigma_{h}^{\frac{3-\nu}{2}}.
\end{align*}
This yields the bound
\begin{equation}\label{eq:finalDelta1}
\Delta_{k}\geq \tilde{\gamma}_{\nu}\lambda_{k}\min\left\{\frac{2\sigma_{h}^{(3-\nu)/2}}{M(\nu-2)},\frac{4-\nu}{(\nu-2)(L+\mu)}\lambda_{k}\right\}\qquad\forall k\geq 0.
\end{equation}
Recall from Proposition \ref{prop:Delta} that $\lim_{k\to\infty}\Delta_{k}=0$ always holds. Consequently, combining \eqref{eq:finalDelta1} with \eqref{eq:norm_equal}, we immediately see
$\lim_{k\to\infty}\lambda_{k}=\lim_{k\to\infty}\beta_{k}=0.$  Now observe that 
\[
\alpha_{k}=\frac{1}{\delta_{k}}\left[1-\left(1+\frac{\delta_{k}}{L+\mu}\frac{4-\nu}{\nu-2}\right)^{\frac{2-\nu}{4-\nu}}\right]=:\frac{Q(\delta_{k})}{\delta_{k}},
\]
and $\delta_{k}=M(\nu/2-1)\lambda_{k}^{\nu-2}\beta_{k}^{3-\nu}$. Thus, 
\begin{equation}\label{eq:delta01}
\lim_{k\to\infty}\delta_{k}=0,
\end{equation}
and $\lim_{k\to\infty}Q(\delta_{k})=0$. By l'H\^{o}pital rule
\begin{equation}\label{eq:alpha1}
\lim_{k\to\infty}\alpha_{k}=\lim_{k\to\infty}\frac{Q(\delta_{k})}{\delta_{k}}=\lim_{k\to\infty}Q'(\delta_{k})=\frac{1}{L+\mu}>0.
\end{equation}
Finally, by definition of the search direction, there exists a sequence of dual variables $(y^{k})_{k\geq 0}\subset\R^{m}$, explicitly defined by \eqref{eq:dual}, for which  
\begin{align*}
\norm{\nabla F_{\mu}(x^{k})-A^{\top}y^{k}}_{x^{k}}^{\ast}=\lambda_{k}\quad\forall k\geq 0.
\end{align*}
We therefore observe first that $\lim_{k\to\infty} \norm{\nabla F_{\mu}(x^{k})-A^{\top}y^{k}}_{x^{k}}^{\ast}=0$, and second 
\begin{align*}
\norm{\nabla F_{\mu}(x^{k})-A^{\top}y^{k}}_{2}\leq \abs{H(x^{k})^{1/2}}\lambda_{k}.
\end{align*}
Since $(x^{k})_{k\geq 0}\subset\scrS_{\mu}(x^{0}),H(x^{k})\succ 0$, and $h\in\bC^{3}(\dom h)$, using \eqref{eq:tau}, we conclude that 
\[
\lim_{k\to\infty}\norm{\nabla F_{\mu}(x^{k})-A^{\top}y^{k}}_{2}\leq\sqrt{\tau_{h}}\lim_{k\to\infty}\lambda_{k}=0.
\]
%%%%%%%%%%%%%%%%%%%%%%%%%%%%%%%%%%%%%%%%%%%%%%%%%%%%%%
\subsection{The case $\nu=3$}

A direct substitution for $\Delta_{k}$ gives us 
\begin{equation}
\label{eq:nu=3_Delta_k}
\Delta_{k}=\frac{\lambda^{2}_{k}}{\frac{M}{2}\lambda_{k}+L+\mu}+\frac{4}{M^{2}}(L+\mu)\left[\frac{\frac{M}{2}\lambda_{k}}{\frac{M}{2}\lambda_{k}+L+\mu}+\ln\left(\frac{L+\mu}{\frac{M}{2}\lambda_{k}+L+\mu}\right)\right]. 
\end{equation}
Denote $t _{k}:= (L+\mu)/(\frac{M}{2}\lambda_{k}),\delta_{k}=\frac{M}{2}\lambda_{k}$. Then 
\[
\alpha_{k} = \frac{2}{M\lambda_k}\frac{1}{1+t_{k}}=\frac{1}{\delta_{k}+L+\mu},
\]
so that 
\[
\frac{\alpha_{k}M\lambda_k}{2}= \frac{1}{1+t_{k}},\text{ and }L+\mu = \frac{M}{2}\lambda_{k}t_{k}.
\]
This implies that 
\begin{align} \label{eq:nu=3_Delta_k_simpl}
\Delta_{k}&=\frac{2\lambda_{k}}{M}\frac{1}{1+t_{k}}+\frac{2\lambda_{k}}{M} t_{k} \left[\frac{1}{1+t_{k}}+\ln\left(\frac{t_{k}}{1+t_{k}}\right)\right], \notag \\
&=\frac{2\lambda_{k}}{M} \left(1+ t_{k} \ln\left(\frac{t_{k}}{1+t_{k}}\right) \right).
\end{align}
Consider the function $\gamma:(0,\infty)\to(0,\infty)$, given by  $\gamma(t):=1+ t \ln\left(\frac{t}{1+t}\right)$. When $t \in (0,1)$, since 
\begin{align*}
\gamma'(t) &= \ln\left(\frac{t}{1+t}\right) + t \frac{1+t}{t}\left(\frac{1}{1+t}-\frac{t}{(1+t)^2}  \right) \\
%& = \ln\left(\frac{t}{1+t}\right) + 1 - \frac{t}{1+t} 
&= \ln\left(1- \frac{1}{1+t}\right)  + \frac{1}{1+t} < 0,
\end{align*}
we conclude that $\gamma(t)$ is decreasing for $t\in(0,1)$. Hence, $\gamma(t) \geq \gamma(1) = 1-\ln 2$, for all $t \in (0,1)$. On the other hand, if $t \geq 1$, 
\begin{align*}
\frac{\dif}{\dif t}\left(\frac{\gamma(t)} {1/t}\right)=\frac{\dif}{\dif t}(t\gamma(t) ) =
%& = 1 + 2t \ln\left(\frac{t}{1+t}\right) + t^2 \frac{1+t}{t}\left(\frac{1}{1+t}-\frac{t}{(1+t)^2}  \right) \\
1 + 2t \ln\left(\frac{t}{1+t}\right)  + \frac{t}{1+t}\geq 0.
\end{align*}
Hence, $t\mapsto \frac{\gamma(t)} {1/t}$ is an increasing function for $t\geq 1$, and thus $\gamma(t) \geq \frac{1-\ln 2}{t}$, for all $t\geq 1$. Summarizing these two cases we see $\Delta_{k}\geq \frac{2\lambda_{k}}{M}\min\{1,1/t_{k}\}(1-\ln(2))$, which after rearranging, can be stated as 
\begin{equation}\label{eq:finalDelta2}
\Delta_{k}\geq \frac{2(1-\ln(2))\lambda_{k}}{M(L+\mu)}\min\left\{L+\mu,\frac{M}{2}\lambda_{k}\right\}\quad\qquad\forall k\geq 0.
\end{equation}
From Proposition \ref{prop:Delta} we know that $\lim_{k\to\infty}\Delta_{k}=0$, and consequently, 
\begin{equation}\label{eq:delta02}
\lim_{k\to\infty}\delta_{k}=\lim_{k\to\infty}\lambda_{k}=0,\text{ as well as } \lim_{k\to\infty}\alpha_{k}=\frac{1}{L+\mu}.
\end{equation}
Eq. \eqref{eq:norm_equal} shows that $\lim_{k\to\infty}\beta_{k}=0$. As in the case $\nu\in(2,3)$, we arrive at the conclusion $\lim_{k\to\infty}\norm{\nabla F_{\mu}(x^{k})-A^{\top}y^{k}}_{2}=0$. 
%%%%%%%%%%%%%%%%%%%%%%%%%%%%%%%%%%%%%%%%%%%%%%%%%%%%%%
\subsection{The case $\nu\in(3,4)$}

Similarly to the case $\nu \in (2,3)$, denote $t_k = 1+\frac{\delta_{k}}{L+\mu}\frac{4-\nu}{\nu-2} \in (1,+\infty)$, $\ca=\frac{4-\nu}{2(3-\nu)}\in(-\infty,0)$, $\cb=\frac{2-\nu}{4-\nu}\in(-\infty,-1)$.  Then the expression for the $\Delta_k$ is the same as for $\nu \in (2,3)$:
\begin{align*}
\Delta_{k}&=\frac{\lambda^{2}_{k}}{\delta_{k}} \left(1+\frac{\ca}{t_k-1}-\frac{\ca t_k^{\cb+1}}{t_k-1}\right). 
\end{align*} 
Let us define a function $\gamma(t)$ such that $\Delta_{k}=\frac{\lambda^{2}_{k}}{\delta_{k}}\gamma(t_k)$. 
Our next goal is to show that, for $t \in [2,+\infty)$, $\gamma(t)$ is below bounded by some positive constant and, for $t \in (1,2]$, $\gamma(t)$ is below bounded by some positive constant multiplied by $t-1$.

\textit{1. $t \in [2,+\infty)$.} We will show that $\gamma'(t) \geq 0$, whence $\gamma(t) \geq \gamma(2)$. Thus, we need to show that
\begin{align*}
0 & \leq \frac{d}{dt} \left( 1+\frac{\ca}{t-1}-\frac{\ca t^{\cb+1}}{t-1}\right) = -\frac{\ca}{(t-1)^2}\underbrace{\left(1-(\cb+1)t^\cb+\cb t^{\cb+1}\right)}_{=:\psi(t)}.
\end{align*}
Since $\ca \leq 0$, to show that $\gamma'(t) \geq 0$ it is enough to show that $\psi(t) \geq 0$. Since $\cb<-1$,
\[
\psi'(t) = \cb(\cb+1)t^\cb - \cb(\cb+1) t^{\cb-1} = \cb(\cb+1)t^{\cb-1}(t-1) \geq 0,
\]
whence, $\psi(t) \geq \psi(2)$, $t \in [2,+\infty)$. It remains to show that $\psi(2) \geq 0$. Let us consider $\psi(2) = 1-(\cb+1)2^\cb+\cb 2^{\cb+1} = 1+\cb 2^\cb-2^\cb$ as a function of $\cb$. For all possible values $\cb\in(-\infty,-1)$ one can check numerically that $\psi(2)\in(0,1)$. Hence, $\psi(t)\geq 0$ for all $t\geq 2$.

\textit{2. $t \in (1,2]$.} We will show that $ \frac{d}{dt} \left(\gamma(t)/(t-1)\right) \leq 0$, whence $\gamma(t) \geq (t-1)\gamma(2)$. Thus, we need to show that
\begin{align*}
0 & \geq \frac{d}{dt} \left( \frac{1}{t-1}+\frac{\ca}{(t-1)^2}-\frac{\ca t^{\cb+1}}{(t-1)^2}\right) \\
&= -\frac{1}{(t-1)^2}-\frac{2\ca}{(t-1)^3}-\frac{\ca(\cb+1)t^\cb}{(t-1)^2}+\frac{2\ca t^{\cb+1}}{(t-1)^3} \\
& = \frac{1}{(t-1)^3} \left( -t + 1 -2 \ca - \ca(\cb+1)t^{\cb+1} + \ca(\cb+1)t^\cb + 2\ca t^{\cb+1} \right) \\
& = \frac{1}{(t-1)^3}\underbrace{\left( -t + 1 -2 \ca + \ca(\cb+1)t^{\cb} - \ca(\cb-1)t^{\cb+1} \right)}_{=:\psi(t)}. 
\end{align*}
Our next step is to show that $\psi(t) \leq 0$. We have
\begin{align*}
\psi'(t) & =-1 + \ca (\cb+1) \cb t^{\ca-1}-\ca(\cb-1)(\cb+1) t^\cb \\
\psi''(t) & =\ca \cb(\cb+1) (\cb-1) t^{\cb-2}-\ca(\cb-1)\cb(\cb+1) t^{\cb-1} \\
& = \ca\cb (\cb+1) (\cb-1) t^{\cb-2} (1-t).
\end{align*}
Using the definition of $\ca,\cb$, and the fact that $\nu \in (3,4)$, we obtain that $\ca(\cb+1) = 1$. Hence, since $t >1$, we obtain that $\psi''(t) \leq 0$. Thus, $\psi'(t) \leq \psi'(1) = 0$, $\psi(t) \leq \psi(1) = 0$, and $\gamma(t)/(t-1) \geq \gamma(2)$.

Combining both cases, we obtain that $\gamma(t) \geq \min\{\gamma(2), (t-1) \gamma(2)\}$. Note that 
\begin{align*}
\gamma(2)\equiv \tilde{\gamma}_{\nu}&:= 1+\frac{4-\nu}{2(3-\nu)}-\frac{4-\nu}{2(3-\nu)}2^{(2-\nu)/(4-\nu)+1}\\
&=1+\frac{1-\exp\left(\frac{2(3-\nu)}{4-\nu}\ln(2)\right)}{\frac{2(3-\nu)}{4-\nu}}
\end{align*}
Via L'H\^{o}spital's rule, once can check that $\lim_{\nu\downarrow 3}\hat{\gamma}(\nu)=1-\ln(2)\in(0,1)$, and $\lim_{\nu\uparrow 4}\hat{\gamma}(\nu)=1$, since $\lim_{\nu\uparrow 4}\frac{2(3-\nu)}{4-\nu}=-\infty$. Consequently, 
\[
\Delta_k\geq \tilde{\gamma}_{\nu} \min\left\{\frac{\lambda^{2}_{k}}{\delta_{k}}, \frac{1}{-\cb}\frac{\lambda^{2}_{k}}{L+\mu}\right\}, \quad\tilde{\gamma}_{\nu}\in(1-\ln(2),1).
\]
By Proposition \ref{prop:Delta}, we know that $\lim_{k\to\infty}\Delta_{k}=0$. Therefore, either $\lambda_{k}\to 0$, or $\frac{\lambda_{k}^{2}}{\delta_{k}}\to 0 $. Suppose there exists $\eps>0$ such that $\lambda_k\geq \eps$ for all $k\geq 0$. Then, 
\begin{align*}
\frac{\lambda_{k}^{2}}{\delta_{k}}=\frac{2}{M(\nu-2)}\lambda_{k}^{4-\nu}\beta_{k}^{\nu-3}\geq \frac{2}{M(\nu-2)}\eps^{4-\nu}\beta_{k}^{\nu-3}.
\end{align*}
Hence, $\beta_{k}\to 0$ must hold. But then eq. \eqref{eq:norm_equal2} implies $\lambda_{k}\to 0$. A contradiction. It follows that $\lambda_{k}\to 0$, and therefore, by \eqref{eq:norm_equal}, $\beta_{k}\to 0$. Using that 
\[
\delta_{k}=M(\nu/2-1)\lambda_{k}^{\nu-2}\beta_{k}^{3-\nu} \stackrel{\eqref{eq:norm_equal2}}{\leq} M(\nu/2-1)\tau_{h}^{\frac{\nu-2}{2}}\beta_{k},
\]
we arrive at the string of inequalities
%\begin{align*}
%\frac{\lambda^{2}_{k}}{\delta_{k}}&=\frac{2}{M(\nu-2)}\lambda_{k}^{4-\nu}\beta_{k}^{\nu-3}\stackrel{\eqref{eq:norm_equal}}{\geq }\frac{2}{M(\nu-2)}\sigma_{h}^{\frac{4-\nu}{2}}\beta_{k}\stackrel{\eqref{eq:norm_equal2}}{\geq} \frac{2}{M(\nu-2)}\left(\frac{\sigma^{4-\nu}_{h}}{\tau_{h}}\right)^{\frac{1}{2}}\lambda_{k}.
%\end{align*}
%\PD{Can we do it simpler since $\nu - 3 \geq 0$?:
\begin{align*}
\frac{\lambda^{2}_{k}}{\delta_{k}}&=\frac{2}{M(\nu-2)}\lambda_{k}^{4-\nu}\beta_{k}^{\nu-3}\stackrel{\eqref{eq:norm_equal2}}{\geq}\frac{2}{M(\nu-2)}\tau_{h}^{-\frac{\nu-3}{2}}\lambda_{k}
\end{align*}
%If this works, we need to change (B.10) and Lemma 9.
%}
Hence, 
\begin{equation}
\lim_{k\to\infty}\delta_{k}=0,\text{ and }\lim_{k\to\infty}\alpha_{k}=\frac{1}{L+\mu}.
\end{equation}
We conclude $\lim_{k\to\infty}\norm{\nabla F_{\mu}(x^{k})-A^{\top}y^{k}}_{2}=0$. Moreover, we obtain the explicit bound
\begin{equation}\label{eq:finalDelta3}
\Delta_{k}\geq \lambda_{k}\tilde{\gamma}_{\nu}\min\left\{\frac{2}{M(\nu-2)}\tau_{h}^{-\frac{\nu-3}{2}},\frac{1}{-\cb}\frac{\lambda_{k}}{L+\mu}\right\}. 
\end{equation}
%%%%%%%%%%%%%%%%%%%%%%%%%%%%%%%%%%%%%%%%%%%%%%%%%%%%%%
\subsection{The case $\nu=4$}
We can compute the per-iteration potential reduction as 
\begin{align*}
\Delta_{k}&=\frac{\lambda_{k}^{2}}{\delta_k}\left[1-\exp\left(-\frac{\delta_{k}}{L+\mu}\right)\right] \notag \\
& -\left(\frac{\lambda_{k}}{\delta_{k}}\right)^{2}(L+\mu)\left(\left[1-\exp\left(-\frac{\delta_{k}}{L+\mu}\right)\right]-\frac{\delta_k}{L+\mu}\exp\left(-\frac{\delta_{k}}{L+\mu}\right)\right). 
\end{align*}
To analyze this expression, denote by $t_k^{-1}:= \frac{\delta_{k}}{L+\mu}$. Then
\begin{align*}
\Delta_{k}&=\frac{\lambda^{2}_{k}}{\delta_{k}} \left( 1-t_k+t_k \exp\left(-\frac{1}{t_k}\right)\right)\geq 0. 
\end{align*} 
Let us define a function $\gamma(t)$ such that $\Delta_{k}=\frac{\lambda^{2}_{k}}{\delta_{k}}\gamma(t_k)$. Our next goal is to show that, for $t \in (0,1]$, $\gamma(t)$ is below bounded by some positive constant and, for $t \geq 1$, $\gamma(t)$ is below bounded by some positive constant divided by $t$.

\textit{1. $t \in (0,1]$}. We will show that $\gamma'(t) \leq 0$, whence $\gamma(t) \geq \gamma(1)$. Indeed, for $t \in (0,1]$,
\begin{align*}
\gamma'(t)&= -1+\exp(-1/t)(1+1/t)\\
&<-1+2\exp(1/t)\\
&\leq-1+2\exp(-1)<0. 
\end{align*}
Thus, we have $\gamma(t) \geq \gamma(1) = \exp(-1).$

\textit{2. $t \in [1,+\infty)$}. We will show that $ \frac{d}{dt} \left(\frac{\gamma(t)}{1/t}\right) \geq 0$, whence $\gamma(t) \geq \frac{\gamma(1)}{t}$.
\begin{equation}
\label{eq:appB_nu=4_proof_1}
\frac{d}{dt} \left( t\left(1-t+t \exp\left(-\frac{1}{t}\right) \right) \right) = \exp\left(-\frac{1}{t}\right) (2t+1) + 1 - 2t .
\end{equation}
Using the Taylor expansion for $\ln(1+x)$ and $\ln(1-x)$ for $x \in (0,0.5]$, we have
\begin{align}
\ln(1+x) - \ln(1-x) &= x - \frac{x^2}{2}+\frac{x^3}{3} + \sum_{k=4}^\infty \frac{(-1)^kx^k}{k} - \left( -x - \frac{x^2}{2}+\frac{x^3}{3} - \sum_{k=4}^\infty \frac{x^k}{k}\right) \notag \\
&= 2x + \frac{2x^3}{3} + \sum_{k=2}^\infty \frac{2x^{2k+1}}{2k+1}  \geq 2x \notag.
\end{align}
Setting $x = \frac{1}{2t}$ for $t\geq 1$, we obtain
\begin{align*}
 & \ln(1+1/(2t)) - \ln(1-1/(2t)) \geq 1/t  \\
\iff & \ln(2t(1+1/(2t))) - \ln(2t(1-1/(2t))) \geq 1/t \\
 \iff & \ln(2t+1) - \ln(2t-1) \geq 1/t\\
\iff & -1/t + \ln(2t+1) \geq \ln(2t-1)\\
\iff & \exp\left(-\frac{1}{t}\right) (2t+1) + 1 - 2t \geq 0.
\end{align*}
which, combined with \eqref{eq:appB_nu=4_proof_1} proves that $ \frac{d}{dt} \left(\frac{\gamma(t)}{1/t}\right) \geq 0$ for $t \geq 1$. Thus, we have that, for $t\geq 1$, $\gamma(t) \geq \frac{\gamma(1)}{t} = \frac{\exp(-1)}{t}$.

Combining two cases $t_k \in (0,1]$ and $t_k \in [1, +\infty)$, we obtain that $\gamma(t) \geq \min\{\gamma(1),\gamma(1)/t\}$ and, since  $t_k^{-1}:= \frac{\delta_{k}}{L+\mu}$
\[
\Delta_k = \frac{\lambda^{2}_{k}}{\delta_{k}}\gamma(t_k) \geq \frac{\lambda^{2}_{k}}{\delta_{k}} \min\{\gamma(1),\gamma(1)/t_k\}= \exp(-1) \min\left\{\frac{\lambda^{2}_{k}}{\delta_{k}},\frac{\lambda^{2}_{k}}{L+\mu}\right\}.
\]

By Proposition \ref{prop:Delta}, we know that $\lim_{k\to\infty}\Delta_{k}=0$. Thus, either $\lambda_{k}\to 0$, or $\frac{\lambda^{2}_{k}}{\delta_{k}}\to 0$. Suppose there exists $\eps>0$ such that $\lambda_{k}\geq\eps>0$ for all $k\geq 0$. Then it must be true $\frac{\lambda^{2}_{k}}{\delta_{k}}\to 0$. Then,  $\frac{\lambda^{2}_{k}}{\delta_{k}}=\frac{\beta_{k}}{M},$ and therefore $\beta_{k}\to 0$ must be true. But then \eqref{eq:norm_equal2} yields the contradiction $\lambda_{k}\to 0$. We are therefore forced to conclude that $\lim_{k\to\infty}\lambda_{k}=0$, and from \eqref{eq:norm_equal} it then follows $\lim_{k\to\infty}\beta_{k}=0$. Furthermore, using \eqref{eq:norm_equal2}, 
\begin{align*}
\delta_{k}=M\frac{\lambda_{k}^{2}}{\beta_{k}}\leq M\tau_{h}\beta_{k},
\end{align*}
so that $\lim_{k\to\infty}\delta_{k}=0$ and $\frac{\lambda_{k}^{2}}{\delta_{k}}\geq \frac{\lambda_{k}}{\sqrt{\tau_{h}}M}.$ This gives the final estimate
\begin{equation}\label{eq:finalDelta4}
\Delta_k \geq \exp(-1)\lambda_{k} \min\left\{\frac{1}{\sqrt{\tau_{h}}M},\frac{\lambda_{k}}{L+\mu}\right\}.
\end{equation}
A simple application of l'H\^{o}pital's rule gives $\lim_{k\to\infty}\alpha_{k}=\frac{1}{L+\mu}$, and $\lim_{k\to\infty}\norm{\nabla F_{\mu}(x^{k})-A^{\top}y^{k}}_{2}=0$.
%%%%%%%%%%%
\begin{Proof}[Proof of Theorem \ref{th:gradient}]
Combining all the results just derived for each generalized self-concordant parameter $\nu\in(2,4]$, we conclude that always $\lim_{k\to\infty}\norm{\nabla F_{\mu}(x^{k})-A^{\top}y^{k}}_{2}=0$. Corollary \ref{cor:interior} shows that $(x^{k})_{k\geq 0}\subseteq\scrS_{\mu}(x^{0})$, which is a compact set by Lemma \ref{lem:S_compact} contained in $\scrX^{\circ}$. Since $\nabla F_{\mu}(x)=\nabla f(x)+\mu\nabla h(x)$ is a continuous function on $\scrX^{\circ}$, we conclude that along every convergent subsequence $(x^{k_{q}})_{q\in\N}$ with limit $\bar{x}\in\scrX^{\circ}$, we have 
\[
\lim_{q\to\infty}A^{\top}y^{k_{q}}=\lim_{q\to\infty}\nabla F_{\mu}(x^{k_{q}})=\nabla F_{\mu}(\bar{x}).
\]
Recall that $y^{k_{q}}=y_{\mu}(x^{k_{q}})$ and the map $x\mapsto y_{\mu}(x)$ is continuous by Lemma \ref{lem:dual}. Denote by $\bar{y}\in\R^{m}$ the corresponding limit of the convergent subsequence $(y^{k_{q}})_{q\in\N}$, we conclude that $\nabla F_{\mu}(\bar{x})=A^{\top}\bar{y}$. Since the convergent subsequence $(x^{k_{q}})_{q\in\N}$ has been chosen arbitrarily, the claim $\omega(x^{0})\subseteq \{x\in\scrX\vert (\exists y\in\R^{m}):\nabla F_{\mu}(x)-A^{\top}y=0\}$ follows. 
\end{Proof}